\documentclass[a4paper,12pt]{article}
\usepackage{a4wide}
\usepackage{theorem}
\usepackage{amsmath}
\usepackage{array}
\usepackage{amssymb}
\usepackage{amsfonts}
\usepackage[french]{babel}
\usepackage{epsf}
\usepackage{epsfig}
\usepackage{graphicx}
\usepackage{float}

\usepackage[french]{babel}
\usepackage{graphicx}
\usepackage{eufrak}

\newtheorem{theorem}{Th\'eor\`eme}[section]
\newtheorem{prop}[theorem]{Proposition}
\newtheorem{lemma}[theorem]{Lemme}
\newtheorem{cor}{Corollaire}
\theoremstyle{definition}

\newtheorem{rem}{Remarque}

\def\N{{\mathbb{N}}}
\def\Z{{\mathbb{Z}}}

\def\R{{\mathbb{R}}}
\def\C{{\mathbb{C}}}
\def\T{{\mathbb{T}}}

\newcommand{\mult}{\mathop{\rm mult}\nolimits}

\newcommand{\goth}[1]{\EuFrak{#1}}
\newcommand{\s}{\mathop{\goth s}\nolimits}

\makeatletter
\newcommand\@makefntextsans[1]{%
    \parindent 0em%
    \noindent%
    \hb@xt@0em{\hss}%
    #1}
\def\footnotetextsans{%
     \@ifnextchar [\@xfootnotenextsans%
       {\@footnotetextsans}}
\def\@xfootnotenextsans[#1]{%
  \begingroup%
     \csname c@\@mpfn\endcsname #1\relax%
  \endgroup%
  \@footnotetextsans}
\long\def\@footnotetextsans#1{\insert\footins{%
    \reset@font\footnotesize%
    \interlinepenalty\interfootnotelinepenalty%
    \splittopskip\footnotesep%
    \splitmaxdepth \dp\strutbox \floatingpenalty \@MM%
    \hsize\columnwidth \@parboxrestore%
    \color@begingroup%
      \@makefntextsans{%
        \rule\z@\footnotesep\ignorespaces#1\@finalstrut\strutbox}
    \color@endgroup}}
\makeatother

\begin{document}

\title{Invariants entiers en g\'eom\'etrie \'enum\'erative r\'eelle}
\author{Jean-Yves Welschinger}
\maketitle

\makeatletter\renewcommand{\@makefnmark}{}\makeatother
\footnotetextsans{Keywords : Enumerative geometry, rational curve, real algebraic variety, holomorphic discs.}
\footnotetextsans{AMS Classification : 53D45, 14N35.}

\begin{abstract}
Je rappelle les divers probl\`emes de g\'eom\'etrie \'enum\'erative r\'eelle desquels j'ai pu extraire des invariants \`a valeurs enti\`eres,
fournissant un pendant r\'eel aux invariants de Gromov-Witten. Je discute l'optimalit\'e des bornes inf\'erieures fournies par ces
invariants ainsi que certaines de leurs propri\'et\'es arithm\'etiques. Je pr\'esente enfin davantage de r\'esultats garantissant la
pr\'esence ou l'absence de disques pseudo-holomorphes \`a bord dans une sous-vari\'et\'e lagrangienne d'une vari\'et\'e symplectique
donn\'ee. 

\end{abstract}

\section*{Introduction}

Le nombre de racines complexes d'un polyn\^ome g\'en\'erique \`a une variable de degr\'e $d$ ne d\'epend pas du choix du polyn\^ome et vaut $d$, tandis que lorsque ce
polyn\^ome est \`a coefficients r\'eels, le nombre de ses racines r\'eelles peut prendre toutes les valeurs de m\^eme parit\'e que $d$ comprises entre $0$ et $d$. Ceci tient au fait que le corps des
nombres complexes est alg\'ebriquement clos au contraire du corps des nombres r\'eels. Bien plus g\'en\'eralement, le nombre de solutions d'un \og syst\`eme de $n$ 
\'equations g\'en\'eriques \fg sur une vari\'et\'e projective complexe lisse de dimension $n$ ne d\'epend que du degr\'e de ces \'equations, alors qu'il d\'epend fortement du choix,
m\^eme g\'en\'erique, de ces \'equations lorsqu'elles sont \`a coefficients r\'eels et consid\'er\'ees sur le lieu r\'eel d'une vari\'et\'e alg\'ebrique r\'eelle. (En fait d'\'equations, il 
conviendrait plut\^ot de parler de sections g\'en\'eriques de $n$ fibr\'es en droites holomorphes disons tr\`es amples). Chaque probl\`eme de g\'eom\'etrie \'enum\'erative r\'eelle 
peut en principe s'interpr\'{e}ter
de cette mani\`ere. La vari\'et\'e projective r\'eelle est l'espace des modules des objets g\'eom\'etriques que l'on veut compter et les \'equations proviennent des conditions
d'incidences que l'on impose \`a ces objets.

Le principal ph\'enom\`ene pr\'esent\'e dans cet article de synth\`ese est le suivant : il est parfois possible de compter ces objets
g\'eom\'etriques r\'eels en fonction d'un signe $\pm$ de mani\`ere \`a extraire un entier ind\'ependant du choix g\'en\'erique des conditions d'incidence. 
Dans le premier paragraphe, nous observons ce ph\'enom\`ene en comptant les courbes $J$-holomorphes rationnelles r\'eelles dans une vari\'et\'e 
symplectique r\'eelle de dimension quatre en fixant leur classe d'homologie et leur imposant de passer par un nombre ad\'equat de points r\'eels ou bien complexes
conjugu\'es. Nous utilisons en effet le langage de la g\'eom\'etrie symplectique pour \'etudier ces probl\`emes \'enum\'eratifs, tenant compte des r\'esultats de M. Gromov \cite{Gro}
selon lesquels le caract\`ere alg\'ebrique des vari\'et\'es ne joue aucun r\^ole dans ces probl\`emes \'enum\'eratifs, seule l'ellipticit\'e de l'op\'erateur de Cauchy-Riemann sous-jacent
intervient. Les entiers que l'on extrait de ce probl\`eme \'enum\'eratif fournissent un invariant par d\'eformation des vari\'et\'es symplectiques r\'eelles de dimension quatre 
$(X, \omega , c_X)$, qui prend la forme d'une fonction $\chi : d \in H_2 (X ; \Z) \mapsto \chi^d [T] \in \Z [T_1 , \dots , T_N]$ o\`u $N$ d\'esigne le nombre de composantes connexes
du lieu r\'eel $\R X$ de la vari\'et\'e. On d\'efinit des invariants analogues pour les vari\'et\'es symplectiques \og fortement semipositives \fg, par exemple positives, dans le 
troisi\`eme paragraphe et en incluant des conditions de tangence \`a une
courbe r\'eelle dans le deuxi\`eme. Ces derniers r\'esultats s'appliquent en particulier \`a un probl\`eme classique de g\'eom\'etrie \'enum\'erative, le comptage des
coniques tangentes \`a cinq coniques g\'en\'eriques donn\'ees. Le nombre de solutions complexes vaut $3264$, un r\'esultat \'etabli par de Joncqui\`eres au milieu du
dix-neuvi\`eme si\`{e}cle. On montre que le nombre de solutions r\'eelles se trouve minor\'e par trente-deux lorsque les coniques r\'eelles bordent cinq disques disjoints
par exemple. En effet, la valeur absolue des invariants entiers que l'on introduit dans ce m\'emoire borne inf\'erieurement le nombre de solutions r\'eelles du probl\`eme
\'enum\'eratif que l'on consid\`ere.

Un deuxi\`eme ph\'enom\`ene appara\^{\i}t dans cet article, l'optimalit\'e de ces bornes inf\'erieures. On montre en effet dans le premier paragraphe
que dans le cas des vari\'et\'es symplectiques r\'eelles de dimension quatre, lorsque le lieu r\'eel poss\`ede une sph\`ere, un tore ou bien, sous des conditions
plus restrictives, un plan projectif r\'eel et lorsqu'au plus un point est choisi r\'eel et dans cette composante, il existe une structure presque complexe g\'en\'erique $J$
pour laquelle le nombre de courbes $J$-holomorphes rationnelles r\'eelles satisfaisant nos conditions d'incidence vaut exactement la valeur absolue de notre invariant,
ceci quelle que soit la classe d'homologie de ces courbes rationnelles. Ce r\'esultat vaut \'egalement pour la quadrique ellipso\"{\i}de de dimension trois, comme \'etabli dans 
le troisi\`eme paragraphe. Cette optimalit\'e est \'etablie \`a l'aide de m\'ethodes issues de la th\'eorie symplectique des champs, m\'ethodes qui nous permettent \'egalement
parfois de calculer le signe de notre invariant, d'\'etablir des congruences satisfaites par ce dernier ainsi que de fournir des formules le calculant dans certains cas, calculs
que l'on m\`ene explicitement en bas degr\'es. Tous ces r\'esultats font l'objet du premier paragraphe de cet article. En utilisant la notion d'involution antibirationnelle
sur une vari\'et\'e symplectique de dimension quatre, on montre de la m\^eme mani\`ere dans le quatri\`eme paragraphe l'existence de disques $J$-holomorphes 
\`a bords dans le tore de Clifford et
satisfaisant des conditions d'incidences ponctuelles. Dans le cas d'une sph\`ere lagrangienne dans une vari\'{e}t\'{e} symplectique n\'egative ou nulle, on montre
au contraire dans ce m\^eme paragraphe, pour tout $E > 0$, l'existence de structures presque-complexes $J$ pour lesquelles aucun disque ou membrane 
$J$-holomorphe d'\'energie inf\'erieure \`a $E$ ne repose sur
cette sph\`ere, un r\'esultat analogue \`a nos r\'esultats d'optimalit\'es puisqu'on atteint ainsi le minimum possible du nombre de disques ou membranes $J$-holomorphes.
Remarquons \`a propos que l'obtention d'invariants entiers ou rationnels \`a partir du comptage des disques $J$-holomophes \`a bords dans une sous-vari\'et\'e
lagrangienne est un probl\`eme classique de g\'eom\'etrie symplectique (et de la th\'eorie des cordes ouvertes en physique th\'eorique) pour lequel peu de solutions existent.
Notre approche en fournit une lorsque la lagrangienne est fix\'ee par une involution antiholomorphe. Remarquons \'egalement que l'absence de disques $J$-holomorphes
pour certaines structures permet de d\'efinir l'homologie de Floer pour des sph\`eres lagrangiennes dans les vari\'et\'es symplectiques \`a premi\`ere classe de Chern nulle,
un autre probl\`eme classique de g\'eom\'etrie symplectique (et de sym\'etrie miroir en physique th\'eorique).

Le pr\'esent article est largement issu de mon m\'emoire d'habilitation \`a diriger des recherches, laquelle fut soutenue \`a l'\'Ecole normale sup\'erieure de Lyon en mars $2008$. \\

{\bf Remerciements :}

Je remercie le Centre national de la recherche scientifique ainsi que l'Agence nationale de la recherche pour leurs soutiens sans lesquels je n'aurais pu r\'ealiser ces travaux.

\section{Invariants \'enum\'eratifs des vari\'et\'es symplectiques r\'eelles de dimension quatre}
\label{sectinv4}

\subsection{D\'efinition des invariants}
\label{subsectinv}

Soit $(X, \omega , c_X)$ une {\it vari\'et\'e symplectique r\'eelle} ferm\'ee de dimension quatre, par quoi on entend une vari\'et\'e symplectique ferm\'ee de dimension quatre
$(X, \omega)$ \'equip\'ee d'une involution $c_X$ satisfaisant la relation $c_X^* \omega = -\omega$. Le lieu fixe $\R X$ de cette involution est suppos\'e ici non-vide, c'est le {\it lieu r\'eel}
de la vari\'et\'e, lequel a la propri\'et\'e d'\^etre lagrangien. Ses composantes connexes sont \'etiquet\'ees $(\R X)_1 , \dots , (\R X)_N$.
On note ${\cal J}_\omega$ l'espace des structures presque-complexes $\omega$-positives de $(X, \omega)$ de classe $C^l$, $l \gg 1$ et $\R {\cal J}_\omega \subset {\cal J}_\omega$
le sous-espace des structures $J$ qui rendent l'involution $c_X$ $J$-antiholomorphe.  Ce sont tous deux des vari\'et\'es de Banach s\'eparables non-vides et contractiles.

Soit $d \in H_2 (X ; \Z)$ une classe d'homologie satisfaisant la relation $(c_X)_* d = -d$ et $J \in \R {\cal J}_\omega$ une structure presque-complexe g\'en\'erique. Les {\it courbes $J$-holomorphes rationnelles r\'eelles} homologues \`a $d$, c'est-\`a-dire les sph\`eres $J$-holomorphes invariantes
par $c_X$ et homologues \`a $d$, forment alors un espace de dimension $c_1(X)d -1$, o\`u $c_1(X)$ d\'esigne la premi\`ere classe de Chern de la vari\'et\'e $(X, \omega)$.
Nous supposons cette dimension positive ou nulle, puisque le cas contraire signifie que l'espace en question est vide, puis faisons
chuter cette dimension \`a z\'ero en imposant quelques contraintes \`a ces courbes, \`a savoir de passer par une collection $\underline{x}$ de $c_1(X)d -1$ points distincts. Ces derniers peuvent
\^etre choisis r\'eels, c'est-\`a-dire fix\'es par $c_X$, ou bien complexes conjugu\'es, c'est-\`a-dire \'echang\'es par $c_X$ ; nous noterons $r_i$ le nombre de points r\'eels choisis dans $(\R X)_i$, 
$i \in \{1, \dots , N \}$, et $r_X$
le nombre de paires de points complexes conjugu\'es, de sorte que $2r_X +\sum_{i=1}^N r_i= c_1(X)d -1$. L'ensemble ${\cal R}_d (\underline{x} , J)$  des courbes $J$-holomorphes rationnelles 
r\'eelles homologues \`a $d$ qui satisfont ces contraintes suppl\'ementaires est fini. Ces courbes sont de plus toutes irr\'eductibles, immerg\'ees et n'ont que des points doubles transverses comme
singularit\'es. Remarquons que le cardinal $R_d (\underline{x} , J) = \# {\cal R}_d (\underline{x} , J)$ d\'epend en g\'en\'eral
des choix auxiliaires de la structure presque complexe et de la configuration de points, essentiellement parce que le corps des r\'eels n'est pas alg\'ebriquement
clos. Nous allons montrer qu'il en devient ind\'ependant lorsque l'on compte ces courbes en fonction d'un signe convenablement choisi.

Soit $C \in {\cal R}_d (\underline{x} , J)$, le nombre total de points doubles de $C$ se calcule par la formule d'adjonction et vaut $\delta = \frac{1}{2} (d^2 - c_1 (X) d +2)$.
Les points doubles r\'eels de $C$ sont de deux natures diff\'erentes. Ils peuvent \^etre l'intersection locale de deux branches r\'eelles ou bien l'intersection locale de deux branches 
complexes conjugu\'ees. Ces points doubles r\'eels sont dits {\it non-isol\'es} dans le premier cas et {\it isol\'es} dans le second

$$\includegraphics{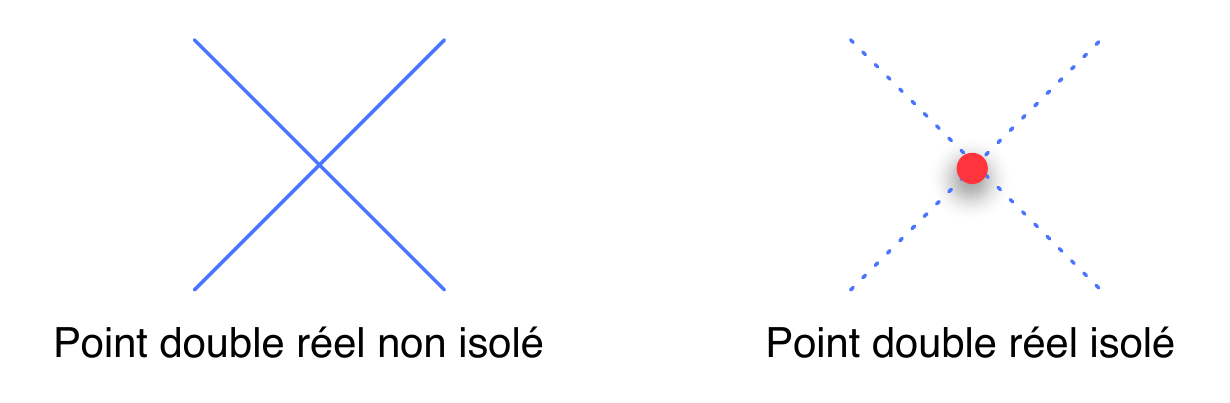}$$

Notons $m(C)$ le nombre de points doubles r\'eels isol\'es de $C$, c'est la {\it masse} de $C$ ; elle est major\'ee par $\delta$. Pour tout entier $m$ compris entre $0$ et $\delta$,
on note $n_d (m)$ le nombre de courbes  $C \in {\cal R}_d (\underline{x} , J)$ de masse $m$. Posons finalement
$$\chi_r^d (\underline{x},J) = \sum_{m=0}^\delta (-1)^m n_d (m),$$
o\`u $r = (r_1 , \dots , r_N)$.
\begin{theorem}[\cite{WelsCRAS1}, \cite{WelsInvent}]
\label{theoInvent}
Soient $(X, \omega , c_X)$ une vari\'et\'e symplectique r\'eelle ferm\'ee de dimension quatre, $N$ le nombre de composantes connexes de son lieu r\'eel et $d \in H_2 (X ; \Z)$ satisfaisant
$c_1 (X)d > 0$. Soient $\underline{x} \subset X$ une configuration r\'eelle de $c_1 (X) d - 1$ points distincts et $r = (r_1 , \dots , r_N)$ le $N$-uplet associ\'e.
L'entier $\chi_r^d (\underline{x},J)$ est ind\'ependant du choix de $\underline{x}$ et du choix g\'en\'erique de $J \in \R {\cal J}_\omega$.
\end{theorem}
Le Th\'eor\`eme \ref{theoInvent} permet de noter cet entier $\chi_r^d$ sans ambigu\"{\i}t\'e. Lorsque $\sum_{i=1}^N r_i$ n'a pas la m\^eme parit\'e que $c_1 (X) d -1$, on pose $\chi_r^d = 0$.
On note alors $\chi^d [T]$ la fonction g\'en\'eratrice $\sum_{|r|=0}^{c_1 (X) d -1} \chi_r^d T^r \in \Z [T_1 , \dots , T_N]$, o\`u $T^r = T_1^{r_1} \dots T_N^{r_N}$ et $|r| = r_1 + \dots + r_N$.
Cette fonction est de m\^eme parit\'e que $c_1 (X) d -1$ et tous ses mon\^omes ne d\'ependent en fait que d'une ind\'etermin\'ee. En effet, la partie r\'eelle d'une sph\`ere holomorphe r\'eelle \'etant connexe, 
l'invariant $\chi_r^d$ est contraint de s'annuler lorsque les points r\'eels de $\underline{x}$ ne sont pas tous choisis dans une m\^eme composante $L$ du lieu r\'eel. 
On adoptera la notation $\chi_r^d (L)$ pour indiquer que les $|r|$ points r\'eels sont choisis dans $L$. On renvoie le lecteur \`a \cite{WelsInvent} pour une \'etude
de la d\'ependance de $\chi_r^d$ en fonction de $r$.

Ainsi, la fonction $\chi : d \in H_2 (X ; \Z) \mapsto \chi^d [T] \in \Z [T_1 , \dots , T_N]$ ne d\'epend que de la vari\'et\'e symplectique r\'eelle ferm\'ee de dimension quatre $(X, \omega , c_X)$ 
et est invariante par d\'eformation de cette derni\`ere. Ceci signifie que si $\omega_t$ est une famille continue de formes symplectiques satisfaisant $c_X^* \omega_t = -\omega_t$, alors
la fonction $\chi$ est la m\^eme pour tous les triplets $(X , \omega_t , c_X)$. Existe-t-il des invariants \'enum\'eratifs analogues \`a ceux qui ressortent du Theor\`eme \ref{theoInvent} en genre quelconque, 
en dimension quelconque et avec des conditions d'incidence quelconques ? Nous n'avons que des d\'ebuts de r\'eponses \`a ces questions. 

\subsection{Bornes inf\'erieures et optimalit\'e}
\label{subsectopt}

Le nombre $R_d (\underline{x} , J)$ de courbes $J$-holomorphes rationnelles 
r\'eelles homologues \`a $d$ qui contiennent l'ensemble $\underline{x}$ de points que l'on s'est donn\'e se retrouve ainsi born\'e inf\'erieurement par la valeur absolue de
l'invariant $\chi_r^d$. Ce nombre est par ailleurs toujours major\'e par le nombre total de courbes $J$-holomorphes rationnelles homologues \`a $d$ et contenant $\underline{x}$, lequel
nombre $N_d$ ne d\'epend ni de $J$ g\'en\'erique, ni de $\underline{x}$ ; c'est un invariant de Gromov-Witten de genre z\'ero de la vari\'et\'e $(X, \omega)$. Ainsi,

\begin{cor}[\cite{WelsInvent}]
\label{corbornes}
Sous les hypoth\`eses du Th\'eor\`eme \ref{theoInvent}, l'encadrement
$$ \vert \chi_r^d  \vert  \leq R_d (\underline{x} , J)   \leq N_d$$
vaut pour tout choix de $\underline{x}$ et tout choix g\'en\'erique de $J \in \R {\cal J}_\omega$. $\square$
\end{cor}

Les bornes inf\'erieures apparaissant dans ce Corollaire \ref{corbornes} se trouvent \^etre parfois optimales. C'est-\`a-dire qu'il est parfois possible d'exhiber une configuration 
$\underline{x}$ et une structure g\'en\'erique $J \in \R {\cal J}_\omega$ telles que toutes les courbes $J$-holomorphes rationnelles r\'eelles compt\'ees par $\chi_r^d$ le sont en fonction
d'un unique et m\^eme signe. Nous pr\'esentons dans ce paragraphe les situations dans lesquelles nous avons \'et\'e en mesure de montrer cette optimalit\'e.

\begin{theorem}[\cite{WelsCRAS2}, \cite{WelsSFT}]
\label{theoopt1}
Soit $(X, \omega , c_X)$ une vari\'et\'e symplectique r\'eelle ferm\'ee de dimension quatre et soit $d \in H_2 (X ; \Z)$ une classe d'homologie satisfaisant $(c_X)_* d = -d$.
Supposons que le lieu r\'eel de cette vari\'et\'e poss\`ede une sph\`ere ou un plan projectif r\'eel $L$. Dans ce dernier cas, supposons que $(X, \omega , c_X)$ est elle-m\^eme
symplectomorphe au plan projectif complexe \'eclat\'e en six boules complexes conjugu\'ees au maximum. Les bornes inf\'erieures apparues dans le Corollaire \ref{corbornes}  
sont sous ces hypoth\`eses optimales d\`es que $0 \leq r \leq 1$.
Le signe de l'invariant $\chi_r^d (L)$ est en outre dans ce cas d\'etermin\'e par l'in\'egalit\'e $(-1)^{\frac{1}{2}(d^2 - c_1(X)d + 2)} \chi^d_r (L) \geq 0$.
\end{theorem}

\begin{rem}
La derni\`ere partie du Th\'eor\`eme \ref{theoopt1} signifie que le signe du coefficient de plus bas degr\'e du polyn\^ome $\chi^d (T)$ introduit au paragraphe \ref{subsectinv} s'interpr\`ete
comme la parit\'e du genre lisse de la classe $d$. Le fait que ce signe puisse \^etre n\'egatif en degr\'es congrus \`a trois ou quatre modulo quatre dans le plan projectif complexe
met en d\'efaut la Conjecture $6$ de \cite{IKS1}. 
\end{rem}

\begin{cor}[\cite{WelsSFT}]
Soit $d$ une classe d'homologie de dimension deux du plan projectif complexe ou de la quadrique ellipso\"{\i}de et $0 \leq r \leq 1$.
Les bornes inf\'erieures (\ref{corbornes}) sont
atteintes pour la structure complexe standard lorsque les points complexes conjugu\'es sont choisis tr\`es proches d'une conique imaginaire pure dans le premier cas
et d'une section hyperplane r\'eelle disjointe de $L$ dans le second. $\square$
\end{cor}

\begin{theorem}[\cite{WelsSFT}]
\label{theoopt2}
Soit $(X, \omega , c_X)$ une vari\'et\'e symplectique r\'eelle ferm\'ee de dimension quatre dont le lieu r\'eel poss\`ede un tore $L$ et soit $d \in H_2 (X ; \Z)$ une classe d'homologie 
satisfaisant $(c_X)_* d = -d$. Les bornes inf\'erieures du Corollaire \ref{corbornes} sont
optimales lorsque $r=1$.
Lorsque le lieu r\'eel est connexe -r\'eduit au tore $L$-, l'invariant $\chi_1^d (L)$ est en outre positif. Dans le cas g\'en\'eral, le signe de l'invariant $\chi_1^d (L)$ est
d\'etermin\'e par l'in\'egalit\'e $(-1)^{\frac{1}{2}(d^2 - c_1(X)d + 2)} \chi^d_1 (L) \geq 0$ lorsque le lieu r\'eel des courbes rationnelles ne s'annule pas dans $H_1 (L ; \Z /2\Z)$,
tandis qu'il est d\'etermin\'e par l'in\'egalit\'e $(-1)^{\frac{1}{2}(d^2 - c_1(X)d + 2)} \chi^d_1 (L) \leq 0$ lorsque ce dernier s'annule.
\end{theorem}

\begin{rem}
Dans le cas particulier de la quadrique hyperbolo\"{\i}de, la positivit\'e de $\chi_1^d (L)$ avait \'et\'e observ\'ee dans \cite{IKS1} par d'autres m\'ethodes.
\end{rem}

De savoir si les bornes sup\'erieures apparues dans le Corollaire \ref{corbornes} sont optimales est un probl\`eme classique de g\'eom\'etrie \'enum\'erative r\'eelle pour lequel on ne sait presque rien.
Le seule chose que je puisse signaler est le crit\`ere suivant.
\begin{cor}[\cite{WelsInvent}]
Sous les hypoth\`eses du Th\'eor\`eme \ref{theoInvent}, supposons que $\chi_r^d $ est positif (resp. n\'egatif). Supposons qu'il existe une configuration r\'eelle de points $\underline{x}$ et une
structure g\'en\'erique $J \in \R {\cal J}_\omega$ telles qu'il existe $\frac{1}{2} (N_d - |\chi_r^d|)$ courbes $J$-holomorphes rationnelles r\'eelles de masses impaires (resp. paires) homologues \`a $d$ 
et passant par $\underline{x}$. Alors, toutes les courbes $J$-holomorphes rationnelles homologues \`a $d$ 
et passant par $\underline{x}$ sont r\'eelles, de sorte que les bornes sup\'erieures du Corollaire \ref{corbornes} sont optimales. $\square$
\end{cor}

 Les bornes inf\'erieures fournies par ces invariants sont-elles optimales en g\'en\'eral ? La question se pose d\'ej\`a dans le
 cas du plan projectif (ou de l'espace projectif de dimension trois, voir le \S \ref{subsectDuke}).

\subsection{Congruences}
\label{sectcong}

\'Etant donn\'ee une classe d'homologie $d \in H_2 (X ; \Z)$ d'une vari\'et\'e symplectique r\'eelle de dimension quatre $(X, \omega , c_X)$, nous noterons 
$g_d = \frac{1}{2} (d^2 - c_1 (X)d + 2)$ le genre lisse de $d$ et $c_d = c_1 (X)d -1$ le degr\'e attendu du polyn\^ome $\chi^d (T)$ d\'efini au \S \ref{subsectinv}.

\begin{theorem}[\cite{WelsSFT}]
\label{theocong1}
Soit $(X, \omega , c_X)$ une vari\'et\'e symplectique r\'eelle ferm\'ee de dimension quatre dont le lieu r\'eel poss\`ede une composante connexe $L$ hom\'eomorphe \`a une sph\`ere.
Soient $d \in H_2 (X ; \Z)$ et $r \in \N$. Lorsque $2r+1 < c_d$, la puissance $2^{\frac{1}{2} (c_d - 2r - 1)}$ divise $\chi^d_r (L)$. 
\end{theorem}

{\bf Exemple :}

Le Th\'eor\`eme \ref{theocong1} s'applique \`a l'ellipso\"{\i}de de dimension deux lorsque $d$ est un multiple positif, disons $\delta > 0$, d'une section plane r\'eelle. 
Dans ce cas, $c_d = 4\delta - 1$ et $g_d = \delta^2 - 2 \delta + 1 = \delta + 1 \mod (2)$. Par cons\'equent,  $2^{2 \delta - r - 1}$ divise $\chi^d_r (L)$ lorsque $r < 2 \delta - 1$.
Nous avons \'egalement montr\'e dans \cite{WelsSFT} que $2^{2 \delta - r}$ divise $\chi^d_r (L)$ lorsque de plus $r = 2 \delta + 1 \mod (4)$ ainsi que la congruence 
$\chi^d_{2 \delta - 3} (L) = 0 \mod (16)$.

\begin{theorem}[\cite{WelsSFT}]
\label{theocong3}
Soit $(X, \omega , c_X)$ une vari\'et\'e symplectomorphe au plan projectif complexe \'eclat\'e en six boules complexes conjugu\'ees au maximum. Soit $d \in H_2 (X ; \Z)$ 
une classe satisfaisant $c_d = c_1 (X)d -1 \geq 0$ et soient $r, r_X$ des entiers naturels satisfaisant la relation $r +  2r_X = c_d$. Lorsque $r+ 1 < r_X$, la puissance
$2^{r_X - r - 1}$ divise $\chi^d_r (L)$. 
\end{theorem}

{\bf Exemple :}

Le Th\'eor\`eme \ref{theocong3} s'applique au plan projectif complexe o\`u $d$ est un multiple positif, disons $\delta > 0$, d'une droite complexe. 
Dans ce cas, $8^{\frac{1}{2}(\delta  - r - 1)}$ divise $\chi^d_r$ lorsque $r+1 < \delta $. Nous avons \'egalement montr\'e dans \cite{WelsSFT} que
$2^{\frac{1}{2}(3\delta  - 3r - 1)}$ divise $\chi^d_r$ lorsque de plus $r = \delta  + 1 \mod (4)$ et $\chi^d_{\delta - 3} = 0 \mod (64)$.

\subsection{Calculs}
\label{subsectcalculs}

L'invariant $\chi^d_r$ qui ressort du Theor\`eme \ref{theoInvent} fut rapidement estim\'e apr\`es que je l'ai introduit. G. Mikhalkin a propos\'e dans  \cite{Mikh}
un algorithme permettant, dans le cas des surfaces toriques r\'eelles, le calcul de cet invariant 
lorsque le nombre $r$ de points choisis r\'eels est maximal. Cet algorithme a \'et\'e plus tard \'etendu par E. Shustin \cite{Shu} pour un choix quelconque de points r\'eels.
Il a \'et\'e utilis\'e par I. Itenberg, V. Kharlamov et E. Shustin \cite{IKS} pour estimer cet invariant, fournissant notamment la minoration $\chi^d_{3d-1} \geq \frac{1}{2} d!$ dans le
cas du plan projectif, le calcul en degr\'e inf\'erieur ou \'egal \`a cinq, puis l'asymptotique $\log \vert \chi^d_{c_1 (X) d-1} \vert \cong \log N_d$ dans le cas des surfaces de Del Pezzo r\'eelles $X$, voir \cite{IKS1}. Ces derniers ont \'egalement
plus r\'ecemment obtenu une formule de type Caporaso-Harris tropicale \cite{IKS2} pour le calcul de $\chi^d_{3d-1}$ dans le plan, apr\`es que A. Gathmann et H. Markwig \cite{Gath} ont obtenus cette formule pour le calcul tropical de $N_d$. 
J. Solomon a \'egalement annonc\'e une formule calculant ces invariants $\chi^d_r$ dans le plan. 
E. Shustin \cite{Shu1} a adapt\'e ces m\'ethodes tropicales pour obtenir des r\'esultats analogues dans le cas de la quadrique ellipso\"{\i}de.
Les m\'ethodes de th\'eorie symplectique des champs que j'ai pour ma part utilis\'e (\cite{WelsCRAS2},  \cite{WelsSFT}) m'ont \'egalement permis d'obtenir des formules de type Caporaso-Harris
mais avec des conditions de tangence imaginaires conjugu\'ees. Les invariants relatifs qui interviennent dans ces formules sont introduits au \S \ref{subsectrelreel}. Je ne rappelle pas ici les formules
g\'en\'erales qui se trouvent dans  \cite{WelsSFT}, mais simplement quelques calculs explicites qui en d\'ecoulent facilement. 

\begin{cor}[\cite{WelsSFT}]
\label{corcalcproj}
Soit $(X, \omega , c_X)$ une vari\'et\'e symplectomorphe au plan projectif complexe. Alors,
$\chi^4 (T) = o(T^2)$, $\chi^5 (T) = 64 + 64T^2 + o(T^3)$,  $\chi^6 (T) = 1024T + 1536T^3 + o(T^4)$,
$\chi^7 (T) = -14336  + 11776T^2 + o(T^3)$ et $\chi^8 (T) = -280576T + o(T^2)$.
\end{cor}

Remarquons que $\chi^3 (T) = 2T^2 + 4 T¬^4 + 8T^8$ ; ce 
calcul de $\chi^d (T)$ en degr\'e trois et les ph\'enom\`enes discut\'es ici s'obtiennent simplement en \'eclatant les neuf points base d'un pinceau de cubiques planes et en 
calculant la caract\'eristique d'Euler du lieu r\'eel de la surface obtenue, 
comme observ\'e par V. Kharlamov \cite{KhDg} d\'ej\`a dans les ann\'ees $90$. Toutefois, m\^eme l'existence d'une quartique rationnelle r\'eelle plane passant par onze 
points r\'eels en position 
g\'en\'erale n'\'etait pas connue avant l'introduction de ces invariants $\chi^d_r$. Les valeurs explicites de $\chi^d_r$ pour $d \leq 9$ et tout $r$ furent entre temps obtenues dans
\cite{ABL} comme cons\'equence d'une formule de type Caporaso-Harris tropicale. Ces r\'esultats mirent en d\'efaut la conjecture de monotonie de \cite{IKS1}, de sorte que la fonction $r \mapsto \chi^d_r$ n'est en g\'en\'eral ni positive, ni monotone. 

\begin{cor}[\cite{WelsSFT}]
\label{corcalc2spher}
Soit $(X, \omega , c_X)$ une vari\'et\'e symplectomorphe \`a la quadrique ellipso\"{\i}de de dimension deux. On note $h$ la classe d'une section plane r\'eelle de bidegr\'e $(1,1)$. Alors,
$\chi^{2h} (T) = 2T^3 + 4T^5 + 6T^7$, $\chi^{3h} (T) = 16T + 16T^2 + o(T^3)$,  $\chi^{4h} (T) = -256T + 320T^3 + o(T^4)$ et $\chi^{5h} (T) = 26880T + o(T^2)$.
\end{cor}

\begin{rem}
\label{remfracrat}
 Cet invariant $\chi^d_r$ peut se d\'efinir purement en termes de fractions rationnelles complexes. Lorsque $r = 4d-1$ par exemple, il compte alg\'ebriquement 
 le nombre de fractions rationnelles $u = P/Q$, $P, Q \in \C [X]$ de degr\'es $d$, modulo reparam\'etrage par les homographies r\'eelles de $PGL_2 (\R)$, telles que l'image 
 $u(\R P^1)$ interpole un ensemble donn\'e g\'en\'erique de $4d-1$ points de la sph\`ere de Riemann. Le signe en fonction duquel il convient de compter ces fractions rationnelles 
 $u$ est pair si $u$ poss\`ede un nombre pair de points critiques dans chaque h\'emisph\`ere $\C P^1 \setminus \R P^1$ et impair sinon. 
 Il serait int\'eressant d'\'etudier cet invariant de la quadrique ellipso\"{\i}de en travaillant uniquement avec des fractions rationnelles.
\end{rem}

Quelle est l'asymptotique de l'invariant $\chi_r^d$, $r \leq 1$, calcul\'e ici ? Nos formules calculent l'invariant
 en fonction d'une somme sur des arbres d\'ecor\'es. Quels sont les arbres qui sont asymptotiquement dominants/n\'egligeables ?
 De plus,  dans le cas du plan projectif par exemple, lorsque $r=3d-1$, notre formule calcule l'invariant comme une somme sur des arbres
 dont certains contribuent positivement et d'autres n\'egativement. Ceci garantit l'existence de structures presque-complexes pour lesquelles
 davantage de courbes rationnelles r\'eelles satisfont nos conditions d'incidences que le nombre impos\'e par l'invariant $\chi^d_{3d-1}$.
 Combien de courbes r\'eelles a-t-on ainsi construit ? Enfin, notre m\'ethode de calcul suivie dans la premi\`ere section s'applique \`a toute vari\'et\'e symplectique 
 de dimension quatre et calcule l'invariant
 $\chi$ en fonction d'invariants de Gromov-Witten de surfaces rationnelles relatifs \`a des courbes de carr\'e $-2$ ou $-4$ lorsque $L$ est une sph\`ere
 ou un plan projectif r\'eel. Que sait-on de ces invariants et qu'en d\'eduire pour l'invariant $\chi$ ? Cette direction de recherche reste \`a d\'evelopper.
 Par ailleurs, j'ignore dans quelles situations exactement il est possible de calculer l'invariant $\chi_r^d$ en fonction d'invariants relatifs imaginaires.
 
\section{Invariants relatifs des vari\'et\'es symplectiques r\'eelles de dimension quatre}

Les invariants $\chi^d_r$ introduits au \S \ref{subsectinv} sont d\'efinis par un  comptage de courbes $J$-holomophes rationnelles r\'eelles soumises \`a des conditions d'incidence ponctuelles. 
Ils forment ainsi un analogue r\'eel aux invariants de Gromov-Witten de genre z\'ero ponctuels. J'ai \'egalement d\'efini de tels invariants en admettant que les courbes soient soumises \`a des conditions
de tangence avec une courbe donn\'ee, dans l'esprit de la th\'eorie des invariants relatifs. Ces conditions de tangence peuvent \^etre r\'eelles ou bien complexe conjugu\'ees. Dans le cas de conditions
r\'eelles, je n'ai pu d\'efinir de tels invariants relatifs qu'en admettant une seule condition de tangence et encore m'a-t-il fallu faire intervenir plusieurs types de courbes singuli\`eres. J'expose ces 
r\'esultats dans le \S \ref{subsectrelreel}. J'ai pu en d\'eduire des bornes inf\'erieures pour le nombre de coniques r\'eelles tangentes \`a cinq coniques donn\'ees, un probl\`eme classique de 
g\'eom\'etrie \'enum\'erative. Dans le cas de conditions de tangence complexes conjugu\'ees, la situation est bien meilleure et de tels invariants peuvent s'obtenir avec les m\^emes m\'ethodes
que celles utilis\'ees au \S \ref{subsectinv}. Je n'ai en fait introduit et utilis\'e ces invariants que dans des cas tr\`es particuliers, en utilisant le langage de la th\'eorie symplectique des champs. 
Ils m'ont \'et\'e utiles pour mener les calculs pr\'esent\'es au \S \ref{subsectcalculs}. J'expose ces r\'esultats dans le \S \ref{subsectrelim}.

\subsection{Invariants relatifs r\'eels}

\subsubsection{D\'efinition des invariants}
\label{subsectrelreel}

Soient $(X, \omega , c_X)$ une vari\'et\'e symplectique r\'eelle de dimension quatre,
$d \in H_2(X ; \Z)$ une classe d'homologie telle que $c_1 (X) d \geq 2$ et $\underline{x}$ une configuration r\'eelle
de $c_1 (X) d -2$ points distincts. Comme au \S \ref{subsectinv}, on note  $\R X_1 , \dots , \R X_N$ les composantes connexes
de $\R X$ et $r_i$ le cardinal de $\underline{x} \cap \R X_i$, $i \in \{ 1, \dots , N \}$. Soit $B \subset \R X$ une surface \`a bord lisse.
En chaque point r\'eel $x_i$ de
$\underline{x}$, on choisit une droite vectorielle $T_i$ dans le plan tangent $T_{x_i} \R X$. Pour toute structure
presque complexe $J \in \R {\cal J}_\omega$ suffisamment g\'en\'erique, 
on d\'efinit l'entier $\Gamma^{d,B}_r (J, \underline{x})$ comme la somme des nombres de courbes $J$-holomorphes rationnelles r\'eelles qui r\'ealisent la 
classe d'homologie $d$, passent par la configuration $\underline{x}$ et qui proviennent des quatre familles suivantes :
\begin{itemize}
\item Les courbes tangentes au bord de $B$, elles sont compt\'ees en fonction de leurs masses et de leur contact int\'erieur ou 
ext\'erieur \`a $B$ au point de tangence.
\item Les courbes non-immerg\'{e}es, qui sont compt\'ees en fonction de leurs masses et de la position du point de rebroussement par rapport \`a $B$.
\item Les courbes poss\'edant une des droites $T_i$ comme tangente, qui sont compt\'ees en fonction de leurs masses
et de la position du point $x_i$ correspondant \`a $T_i$ par rapport \`a $B$.
\item Les courbes r\'eductibles, qui sont compt\'ees en fonction de leurs masses
et d'une multiplicit\'e qui est le nombre de points r\'eels d'intersection
entre les deux composantes irr\'eductibles de la courbe, chacun de ces points devant \^etre compt\'e positivement ou 
n\'egativement selon qu'il est int\'erieur ou ext\'erieur \`a $B$. \\
\end{itemize}

Ainsi, en notant respectivement ${\cal T}an^d_B (J , \underline{x})$, ${\cal C}usp^d (J , \underline{x})$,  ${\cal T}an^d (J , \underline{x})$ et ${\cal R}ed^d (J , \underline{x})$
ces quatre ensembles finis de courbes $J$-holomorphes, l'entier $\Gamma^{d,B}_r (J , \underline{x}) $ s'\'ecrit

$$ \sum_{C \in \cup {\cal T}an^d_L (J , \underline{x}) \cup {\cal T}an^d (J , \underline{x}) \cup {\cal C}usp^d (J , \underline{x})}
(-1)^{m(C)} \langle C,B \rangle - \sum_{C \in {\cal R}ed^d (J , \underline{x})} (-1)^{m(C)} \mult_B (C).$$

Dans cette somme, {\it l'indice de contact} $\langle C,B \rangle$ vaut $-1$ (resp. $+1$) si $C \in {\cal T}an^d_L (J , \underline{x})$
et $\R C$ se trouve localement incluse dans (resp. en dehors de) $B$ au voisinage du point de tangence $y$ avec $\partial B$.
Si $C \in {\cal C}usp^d (J , \underline{x})$ (resp. $C \in {\cal T}an^d (J , \underline{x})$), le point de rebroussement (resp. la droite $T_i$, $i \in I$)
est unique et l'indice de contact $\langle C,B \rangle$ vaut $-1$ si ce point se situe en-dehors de $B$ et $+1$ sinon.
Si $C$ est r\'eductible, elle n'a que deux composantes irr\'eductibles $C_1$, $C_2$, toutes deux r\'eelles et 
$$\mult_B (C) = \sum_{y \in \R C_1 \cap \R C_2}  \langle y,B \rangle,$$
o\`u $ \langle y,B \rangle$ vaut $-1$ lorsque $y$ est ext\'erieur \`a $B$ et $+1$ s'il est int\'erieur.

\begin{theorem}[\cite{WelsGAFA}]
\label{theoGAFA}
Soient $(X, \omega , c_X)$ une vari\'et\'e symplectique r\'eelle ferm\'ee de dimension quatre et $B \subset \R X$ une surface \`a bord lisse.
Soient $N$ le nombre de composantes connexes de $\R X$ et $d \in H_2 (X ; \Z)$ satisfaisant $c_1 (X)d > 1$, $c_1 (X)d \neq 4$. 
Soient $\underline{x} \subset X \setminus \partial B$ une configuration r\'eelle de $c_1 (X) d - 2$ points distincts et $r = (r_1 , \dots , r_N)$ le $N$-uplet associ\'e,
suppos\'e non nul.
L'entier $\Gamma^{d,B}_r (J , \underline{x}) $  est ind\'ependant du choix de $\underline{x}$ et du choix g\'en\'erique de $J \in \R {\cal J}_\omega$.
\end{theorem}
Le Th\'eor\`eme \ref{theoGAFA} permet sans ambigu\"{\i}t\'e de noter $\Gamma^{d,B}_r$ cet entier. Lorsque $\sum_{i=1}^N r_i$ n'a pas la m\^eme parit\'e que $c_1 (X) d$, 
on pose $\Gamma^{d,B}_r = 0$. Comme au \S \ref{subsectinv}, on note alors $\Gamma^{d,B} [T]$ la fonction g\'en\'eratrice 
$\sum_{|r|=0}^{c_1 (X) d -2} \Gamma^{d,B}_r T^r \in \Z [T_1 , \dots , T_N]$.
Cette fonction est de m\^eme parit\'e que $c_1 (X) d$ et tous ses mon\^omes ne d\'ependent en fait que d'une ind\'etermin\'ee. 

Ainsi, la fonction $\Gamma^B : d \in H_2 (X ; \Z) \mapsto \Gamma^{d,B} [T] \in \Z [T_1 , \dots , T_N]$ ne d\'epend que du quadruplet $(X, \omega , c_X, B)$. Elle est en outre
invariante par d\'eformation de ce quadruplet au sens o\`u si $\omega_t$ est une famille continue de formes symplectiques satisfaisant $c_X^* \omega_t = -\omega_t$
et $B_t \subset \R X$ une isotopie de surfaces compactes, alors cette fonction est la m\^eme pour tout $(X, \omega_t , c_X, B_t)$. \\

Remarquons qu'en particulier $\Gamma^{d,B}_r (J , \underline{x}) $ ne d\'epend pas de la position relative de $\underline{x}$ par rapport \`a $B$, que 
$\Gamma^{d,B}_r = - \Gamma^{d,\R X \setminus B}_r$ et que le cas particulier o\`u $B$ est vide est admissible et fournit un invariant que l'on a pr\'ealablement
introduit dans \cite{WelsSMF}. Montrer l'invariance de $\Gamma^{d,\emptyset }_r (J , \underline{x}) $ se trouve \^etre une \'etape importante dans la d\'emonstration de
l'invariance de $\Gamma^{d,B}_r (J , \underline{x}) $.

\begin{theorem}[\cite{WelsGAFA}]
\label{theorelation}
Sous les hypoth\`eses du Th\'eor\`eme \ref{theoGAFA}, si $B$ est un disque, 
$2\chi^d_{r+1} = \Gamma^{d,B}_r - \Gamma^{d,\emptyset }_r $. Si de plus,  $(X, \omega , c_X)$ 
est symplectomorphe au plan projectif complexe, $\Gamma^{d,B}_r = -\Gamma^{d,\emptyset }_r $, tandis que
si elle est symplectomorphe \`a la quadrique hyperbolo\"{\i}de de dimension deux, $\Gamma^{d,B}_r = 2\chi^d_{r+1} - \Gamma^{d,\emptyset }_r $.
\end{theorem}

\begin{cor}[\cite{WelsGAFA}]
Sous les hypoth\`eses du Th\'eor\`eme \ref{theorelation}, 
$\chi^d_{r+1} = - \Gamma^{d,\emptyset }_r  =  \Gamma^{d,B}_r$ dans le cas du plan projectif complexe et
$\Gamma^{d,B}_r = 2 \chi^d_{r+1}$, $\Gamma^{d,\emptyset }_r  = 0$ dans le cas de la quadrique hyperbolo\"{\i}de de dimension deux. $\square$
\end{cor}

\subsubsection{Sur les $3264$ coniques tangentes \`a cinq coniques g\'en\'eriques}

Il est possible d'\'etendre les r\'esultats du \S \ref{subsectrelreel}  \`a davantage de conditions de tangence avec le bord de $B$, au moins dans le cas de coniques. 
J'ai illustr\'e ce ph\'enom\`ene en m'int\'eressant au probl\`eme
ancien du comptage du nombre de coniques tangentes \`a cinq coniques g\'en\'eriques donn\'ees. Le nombre de solutions
complexes est $3264$ comme l'a d\'emontr\'e de Joncqui\`eres en $1859$ mais le nombre de solutions r\'eelles d\'epend
du choix des cinq coniques g\'en\'eriques. Soient $B_1 , \dots , B_5$ cinq disques plong\'es dans $\R P^2$
de sorte que leurs bords soient transverses deux \`a deux, et $J \in \R {\cal J}_\omega$. Notons $\Gamma^B$ le nombre de coniques
$J$-holomorphes r\'eelles qui sont soit :
\begin{itemize}
\item irr\'{e}ductibles,  tangentes \`a $B_1 , \dots , B_5$, et compt\'ees positivement si elles sont tangentes int\'erieurement 
\`a $B_i$ pour un nombre pair de $i \in \{ 1, \dots , 5\}$, n\'egativement sinon.
\item r\'eductibles, tangentes \`a quatre des cinq disques $B_1 , \dots , B_5$, et chacune compt\'ee en fonction de la parit\'e
du nombre de disques en lesquelles elle est tangente int\'erieurement et de la position de son unique point singulier par rapport 
au cinqui\`eme disque en lequel elle n'est pas tangente.\\
\end{itemize}

Ainsi, en notant respectivement ${\cal C}on (J)$ et ${\cal C}on_{red} (J)$ ces deux ensembles finis de coniques $J$-holomorphes, on obtient
$$\Gamma^B (J) = \sum_{C \in {\cal C}on (J)} \langle C,B \rangle - \sum_{C \in {\cal C}on_{red} (J)} \langle C,B \rangle \mult_B (C) \in \Z,$$

o\`u lorsque  $C \in {\cal C}on (J)$, {\it l'indice de contact} $\langle C,B \rangle$ vaut $\Pi_{i=1}^5 \langle C,B_i \rangle$ ; tandis que lorsque $C \in {\cal C}on_{red} (J)$ et $i_1 , \dots ,
i_4 \in \{ 1, \dots , 5 \}$ sont les entiers tels que  $C$ soit tangent aux bords de  $B_{i_1} , \dots , B_{i_4}$, $\langle C,B \rangle = \Pi_{j=1}^4 \langle C,B_{i_j}  \rangle$ 
et $\mult_B (C) = +1$ si le point singulier de $C$ appartient \`a $B_{i_5}$ et $-1$ sinon.

\begin{theorem}[\cite{WelsGAFA}]
\label{theocon}
1) Cet entier $\Gamma^B (J)$ ne d\'epend pas du choix g\'en\'erique de la structure presque-complexe $J \in \R {\cal J}_\omega$ et est invariant par
isotopie de $B = B_1 \cup \dots \cup B_5$.

2) Si $B_1 , \dots , B_5$ sont cinq disques disjoints, alors $\Gamma^B = 272$. Il en est de m\^eme si $B_1 , \dots , B_5$
sont proches de cinq droites doubles g\'en\'eriques. 
\end{theorem}

Un disque est dit proche d'une droite double d'\'equation $y^2 = 0$ s'il a une \'equation de la forme
$\{ y^2 \leq \epsilon^2 x^2 - \delta \}$ pour $\epsilon$ et $\delta$ petits.

$$\includegraphics{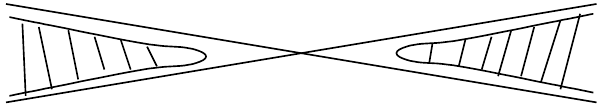}$$

\begin{cor}[\cite{WelsGAFA}]
\label{coninf}
Si $C_1 , \dots , C_5$ sont cinq coniques dont la classe d'isotopie est donn\'ee par la deuxi\`eme partie du Th\'eor\`eme
\ref{theocon}, alors le nombre de coniques r\'eelles qui leur sont tangentes est minor\'e par $32$ ind\'ependamment
du choix de $C_1 , \dots , C_5$ dans la classe d'isotopie.
\end{cor}

{\bf D\'emonstration :}

Le nombre de droites tangentes \`a deux coniques g\'en\'eriques vaut quatre, elles sont cod\'ees par les points d'intersection entre les deux
coniques duales. Le nombre de coniques tangentes \`a quatre des cinq coniques $C_1 , \dots , C_5$ se trouve donc major\'e par $240 = 5*3*4*4$, 
de sorte que le r\'esultat d\'ecoule de la d\'efinition de $\Gamma^B$ et de la deuxi\`eme partie du Th\'eor\`eme \ref{theocon}. $\square$\\

Remarquons  que le fait qu'il existe une configuration de cinq coniques r\'eelles pour lesquelles les $3264$ coniques tangentes \`a ces cinq coniques sont toutes
r\'eelles a \'et\'e \'etabli par F. Ronga, A. Tognoli et T. Vust \cite{Ronga}. Le th\'eor\`eme \ref{theoGAFA} montre la difficult\'e \`a d\'efinir des invariants relatifs avec conditions
de tangence r\'eelles. Dans le \og monde tropical \fg, la situation est parfois bien meilleure, voir \cite{IKS2}. Il en est de m\^eme avec des conditions de tangence complexes
conjugu\'ees, voir le \S \ref{subsectrelim}.

\subsection{Invariants relatifs imaginaires}
\label{subsectrelim}

Soit $L$ une sph\`ere, un tore ou un espace projectif r\'eel de dimension $n=2$ ou $3$. Le fibr\'e cotangent
de $L$ est \'equip\'e de sa forme de Liouville $\lambda$ et de l'involution $c_L$ d\'efinie par $(q,p) \in T^* L \mapsto (q,-p) \in T^* L$. 
Cette derni\`ere satisfait $c_L^* \lambda = - \lambda$ de sorte que $(T^* L , d\lambda , c_L)$ est une vari\'et\'e symplectique r\'eelle.
Soit $g$ une m\'etrique \`a courbure constante sur $L$, $U^* L$  l'ensemble des couples $(q,p) \in T^* L$ tels que $g(p,p) \leq 1$ et $S^* L$ 
le bord de $U^* L$. La restriction de $\lambda$ \`a $S^* L$ est une forme de contact et l'on note $R_\lambda$
le champ de Reeb associ\'e. Le flot engendr\'e par $R_\lambda$ n'est autre que le flot g\'eod\'esique.
Notons ${\cal J}_\lambda$ l'espace des structures presque-complexes positives pour $d \lambda$ et asymptotiquement cylindriques
sur une structure $CR$ de $S^* L$.   Plus pr\'ecis\'ement, le champ radial de $ T^* L$ identifie le compl\'ementaire de la section nulle
avec la symplectisation $(\mathbb{R} \times S^* L , d(e^\rho \lambda))$ de $(S^* L , \lambda)$. On note ${\cal J}_\lambda$ 
l'espace des structures presque-complexes $J$ positives pour $d \lambda$, de classe $C^l$, $l \gg 1$,  qui satisfont
$J(\frac{\partial}{\partial \rho}) = R_\lambda$ et pr\'eservent le noyau de $\lambda$ pour  $\rho \gg 1$ et qui enfin sont invariantes par translation par $\rho$ au-del\`a
d'un certain rang $\rho_0$. Nous notons alors $\mathbb{R} {\cal J}_\lambda \subset {\cal J}_\lambda$ le sous-espace des structures
presque-complexes pour lesquelles $c_L$ est $J$-antiholomorphe. Ces espaces ${\cal J}_\lambda$ et $\Bbb{R} {\cal J}_\lambda$ 
sont tous deux des vari\'et\'es de Banach s\'eparables non-vides et contractiles.
Nous allons compter les courbes $J$-holomorphes rationnelles r\'eelles point\'ees d'\'energie de Hofer finie proprement immerg\'ees dans $T^* L$ en fonction d'un signe $\pm 1$
de fa\c{c}on \`a obtenir un invariant associ\'e \`a $T^* L$. Rappelons que d'apr\`es le Th\'eor\`eme $1.2$ de \cite{HWZ1} et d'apr\`es \cite{Bour}, ces courbes rationnelles point\'ees
convergent en leurs pointes vers des orbites de Reeb parcourues un nombre entier de fois, que l'on appelle multiplicit\'e.
La dimension de l'espace des modules de telles courbes d\'epend du nombre de pointes
et des multiplicit\'es associ\'ees. Afin d'obtenir un nombre fini de courbes, nous allons soumettre ces courbes \`a quelques contraintes, soit en les
for\c{c}ant \`a converger vers des orbites de Reeb prescrites, soit en les for\c{c}ant \`a passer par des points de $L$ ou des paires de points
complexes conjugu\'ees de $T^* L \setminus L$.

Soit $e_i$, $i \geq 1$, la suite d'entiers partout nulle sauf au $i$-\`eme rang o\`u elle vaut un. Soient $\alpha = \sum_{i \in \Bbb{N}^*} \alpha_i e_i$ 
et $\beta = \sum_{i \in \Bbb{N}^*} \beta_i e_i$ deux suites d'entiers positifs qui s'annulent \`a partir d'un certain rang. Ces deux suites
codent respectivement le nombre de paires d'orbites de Reeb complexes conjugu\'ees limites prescrites et non prescrites de nos courbes, avec leur multiplicit\'es $i \in \Bbb{N}^*$.
Le nombre de pointes de nos courbes vaut donc $2v = 2 \sum_{i \in \Bbb{N}^*} (\alpha_i + \beta_i )$ et nous choisissons un ensemble 
$\Gamma$ de $\sum_{i \in \Bbb{N}^*} \alpha_i $ g\'eod\'esiques ferm\'ees disjointes de $L$ pour prescrire nos paires d'orbites de Reeb limites.
\`A pr\'esent, afin de fixer nos contraintes ponctuelles, soient  $r \in \Bbb{N}$ et $x_1 , \dots , x_r$  des points distincts de $L$. 
De m\^eme, soient $r_L \in \Bbb{N}$ et $\xi_1, \overline{\xi}_1 , \dots ,
\xi_{r_L}, \overline{\xi}_{r_L}$ des paires distinctes de points complexes conjugu\'es de $T^* L \setminus L$, c'est-\`a-dire satisfaisant $c_L (\xi_i) = \overline{\xi}_i$.
Nous supposons que
\begin{equation}
\label{dimsphere}
(n-1)r + 2(n-1)r_L + 2(n-1) \# \Gamma = 2 v + \epsilon (n-1) \sum_{i \in \Bbb{N}^*} i (\alpha_i + \beta_i ) + n - 3,
\end{equation}
o\`u $\epsilon = 2$ si $L$ est hom\'eomorphe \`a une sph\`ere et  $\epsilon = 1$ si $L$ est hom\'eomorphe \`a un espace projectif r\'eel, tandis que nous supposons
\begin{equation}
\label{dimtore}
(n-1)r + 2(n-1)r_L  = 2 v + n - 3 \text{ et } \alpha = 0
\end{equation}
si $L$ est hom\'eomorphe \`a un tore.

Alors, lorsque la structure presque-complexe $J \in \Bbb{R} {\cal J}_\lambda$ est g\'en\'erique, il n'y a qu'un nombre fini
de courbes $J$-holomorphes rationnelles r\'eelles d'\'energie de Hofer finie, proprement immerg\'ees dans $T^* L$
et ayant $2v$ pointes qui passent par $\underline{x}$, par chaque paire $\{ \xi_i ,
\overline{\xi}_i \}$ et qui convergent vers les orbites de Reeb relevant les \'el\'ements de $\Gamma$  ainsi que vers
$\beta_j$ autres paires d'orbites, $j \in \Bbb{N}^*$, chacune avec multiplicit\'e $j$ ou de classe d'homologie donn\'ee si $L$ est un tore.
En effet, si $L$ est un tore, il y a une infinit\'e de g\'eod\'esiques ferm\'ees primitives non homologues et la dimension (\ref{dimtore}) ne d\'epend pas du choix
des classes d'homologies de sorte qu'il y a une infinit\'e d'espaces de modules ayant la m\^eme dimension. Pour garantir la finitude, nous imposons les classes
 d'homologies des orbites de Reeb limites. 
Notons ${\cal R} (\alpha , \beta , \Gamma , \underline{x} , \underline{\xi} , J)$ cet ensemble fini de courbes,
la g\'en\'ericit\'e de $J$ garantit qu'elles sont toutes immerg\'ees. 
Si $L$ est de dimension deux, on pose 
$$F_{(r, r_L)} (\alpha , \beta , \Gamma , \underline{x} , \underline{\xi} , J) = \sum_{C \in {\cal R} (\alpha , \beta , \Gamma , \underline{x} , \underline{\xi} , J)} (-1)^{m(C)} \in \Bbb{Z}.$$
Si $L$ est de dimension trois, on l'\'equipe d'une structure spin. Ceci permet d'associer un \'etat spinoriel $\text{sp} (C)$ \`a chaque courbe 
$C \in {\cal R} (\alpha , \beta , \Gamma , \underline{x} , \underline{\xi} , J)$ comme expliqu\'e au \S \ref{subsectDuke} et on pose
$$F_{(r, r_L)} (\alpha , \beta , \Gamma , \underline{x} , \underline{\xi} , J) = \sum_{C \in {\cal R} (\alpha , \beta , \Gamma , \underline{x} , \underline{\xi} , J)} \text{sp} (C) \in \Bbb{Z}.$$

\begin{theorem}[\cite{WelsSFT}]
\label{theoinv}
Soit $L$ une sph\`ere, un tore ou un espace projectif r\'eel de dimension $n=2$ ou $3$ muni d'une m\'etrique \`a courbure constante. Soient
$\alpha , \beta$ deux suites d'entiers positifs qui s'annulent \`a partir d'un certain rang. On choisit comme ci-dessus
un ensemble $\Gamma$ de g\'eod\'esiques ferm\'ees et des ensembles $\underline{x}$, $\underline{\xi}$
de $r$ et $r_L$ points dans $L$ et $T^* L \setminus L$ respectivement de sorte que ces nombres satisfassent (\ref{dimtore}) dans le cas du tore et (\ref{dimsphere}) sinon. 
Lorsque $n=3$, on suppose $r \neq 0$ et lorsque de plus $L \in \{ S^3 , \R P^3 \}$, on suppose que $J$ est invariante par le flot de Reeb pour $\rho \gg 1$. Alors, l'entier $F_{(r, r_L)} (\alpha , \beta , \Gamma , \underline{x} , \underline{\xi} , J)$ d\'efini ci-dessus ne d\'epend ni du choix 
des contraintes $\Gamma , \underline{x} , \underline{\xi}$, ni du choix g\'en\'erique de la structure presque-complexe $J \in \Bbb{R} {\cal J}_\lambda$.
\end{theorem}

L'entier $F_{(r, r_L)} (\alpha , \beta , \Gamma , \underline{x} , \underline{\xi} , J)$ \'etant ind\'ependant de $\Gamma , \underline{x} , \underline{\xi} , J$ d'apr\`es le Th\'eor\`eme \ref{theoinv},
nous le noterons $F_{(r, r_L)} (\alpha , \beta)$.
Afin d'all\'eger encore cette notation, nous noterons cet entier $F (\alpha , \beta)$ lorsque $r_L = 0$,  puisque la valeur de $r$ est alors d\'efinie sans ambigu\"{\i}t\'e par les calculs
de dimensions (\ref{dimsphere}) et  (\ref{dimtore}). Les Lemmes \ref{lemmecalc1}, \ref{lemmecalc2} et \ref{lemmecalc3} fournissent quelques calculs que l'on a pu mener. Les r\'esultats du 
\S  \ref{subsectcalculs} reposent sur ces calculs. 

\begin{lemma}[\cite{WelsSFT}]
\label{lemmecalc1}
Si $L$ est hom\'eomorphe \`a une sph\`ere de dimension deux et $r_L = 0$, on a
$F (e_1 , 0) = F (0 , e_1) = 1$, $F (e_2 , 0) = 2$, $F (0 , e_2) = 8$, $F (2e_1 , 0) = 2$,
$F (e_1 , e_1) = 4$ et $F (0 , 2e_1) = 6$.
\end{lemma}

\begin{lemma}[\cite{WelsSFT}]
\label{lemmecalc2}
Si $L$ est hom\'eomorphe \`a un plan projectif r\'eel  et $r_L = 0$, on a
$F (e_1 , 0) = F (0 , e_1) = F (e_2 , 0) = F (2e_1 , 0) = F (e_1 , e_1)  = F (0 , 2e_1) = 1$
et $F (0 , e_2 ) = 4$.
\end{lemma}

\begin{lemma}[\cite{WelsSFT}]
\label{lemmecalc3}
Si $L$ est hom\'eomorphe \`a un plan projectif r\'eel  et $r_L = 0$, on a 
$F (e_3 , 0) = 2$, $F (0, e_3) = 12$, $F (e_1 + e_2 , 0) = 2$, $F (e_1 , e_2) = 8$,
$F (e_2 , e_1) = 4$, $F (0 , e_1 + e_2) = 24$, $F (3e_1 , 0) = 2$, $F (2e_1 , e_1) = 4$,
 $F (e_1 , 2e_1) = 6$ et $F (0 , 3e_1) = 8$.
\end{lemma}

Toutefois, la valeur de l'invariant $F$ qui ressort du Th\'eor\`eme \ref{theoinv} n'est pas connue en g\'en\'eral. Il serait int\'eressant de d\'evelopper des m\'ethodes permettant son
calcul. 

\section{Invariants en dimension six}

Nous exposons dans ce paragraphe les r\'esultats analogues \`a ceux pr\'esent\'es au \S \ref{sectinv4} que l'on a pu \'etablir en dimension six.

\subsection{D\'efinition des invariants dans les vari\'et\'es alg\'ebriques r\'eelles convexes}
\label{subsectDuke}

Rappelons qu'une vari\'et\'e projective lisse est dite convexe lorsque le groupe $H^1 (\C P^1 ; u^* TX)$
s'annule pour tout morphisme $u : \C P^1 \to X$. Les principaux exemples que je connaisse sont les espaces homog\`enes projectifs,
citons les produits d'espaces projectifs, la quadrique de $\C P^4$ ou encore la vari\'et\'e des drapeaux de $\C^3$.
Il est \`a nouveau possible de d\'efinir un invariant en comptant alg\'ebriquement le nombre de courbes rationnelles r\'eelles qui
r\'ealisent une classe d'homologie $d$ donn\'ee et passent par une configuration r\'eelle de points $\underline{x}$
de cardinal $\frac{1}{2} c_1 (X) d$, o\`u $c_1 (X)$ d\'esigne la premi\`ere classe de Chern de la vari\'et\'e et
$c_1 (X) d$ est suppos\'e pair. Toutefois, le signe
$\pm 1$ en fonction duquel il convient de compter les courbes rationnelles r\'eelles est plus d\'elicat \`a d\'efinir. Le lieu r\'eel
$\R X = \text{fixe} (c_X)$ de $X$ est une vari\'et\'e lisse de dimension r\'eelle trois que l'on suppose orientable pour simplifier. 
Munissons-la d'une orientation ainsi que d'une m\'etrique
riemannienne auxiliaire. Son $SO_3 (\R)$-fibr\'e principal des rep\`eres orthonorm\'es directs s'\'etend alors en un $Spin_3$-fibr\'e 
principal. En effet, l'obstruction \`a l'existence d'une telle extension est donn\'ee en g\'en\'eral par la classe
caract\'eristique $w_2 (\R X)$ et cette obstruction s'annule en dimension trois comme il d\'ecoule
de la formule de Wu. 
Lorsque la configuration r\'eelle
de points est suffisamment g\'en\'erique et poss\`ede au moins un point r\'eel, d'une part les courbes rationnelles
r\'eelles qui passent par $\underline{x}$ et r\'ealisent $d$ sont toutes immerg\'ees (m\^eme lisses en g\'en\'eral) et de 
partie r\'eelle non vide, et d'autre part elle sont \'equilibr\'ees. Ce dernier point signifie que le fibr\'e normal de ces 
courbes se d\'ecompose sur $\C$ comme la somme directe de deux fibr\'es en droite isomorphes $L$ et $M$, fibr\'es qui de 
plus peuvent \^etre choisis r\'eels. Notons ${\cal R}_d (\underline{x})$ cet ensemble fini de courbes rationnelles r\'eelles.
Chaque lieu r\'eel de ces courbes fournit donc un n\oe{}ud immerg\'e dans la vari\'et\'e de
dimension trois $\R X$, et ce n\oe{}ud est de plus canoniquement \'equip\'e d'un rep\`ere mobile ou plut\^ot d'axes mobiles
donn\'es par la tangente au n\oe{}ud et les lieux r\'eels des fibr\'es en droites $L$ et $M$ (en fait, seule la classe 
d'homotopie de ces axes mobiles est canoniquement d\'efinie, puisque la d\'ecomposition du fibr\'e normal en somme $L
\oplus M$ n'est pas uniquement d\'efinie, mais c'est amplement suffisant pour nos besoins). Lorsque les lieux r\'eels de $L$ et $M$
sont orientables, c'est-\`a-dire lorsque ces fibr\'es sont de degr\'e pair, ces axes mobiles peuvent \^etre enrichis
de rep\`eres orthonorm\'es. Ainsi, les n\oe{}uds d\'efinis par les courbes rationnelles r\'eelles sont tous \'equip\'es 
de rep\`eres
orthonorm\'es mobiles qui fournissent des lacets dans le $SO_3 (\R)$-fibr\'e principal des rep\`eres orthonorm\'es, lacets qui
rel\`event les n\oe{}uds de $\R X$. Vient alors l'alternative suivante pour chaque lacet : soit ce lacet du $SO_3 (\R)$-fibr\'e
principal des rep\`eres se rel\`eve en un lacet du $Spin_3$-fibr\'e principal donn\'e par la structure $Spin$
$\s$, soit non.
Ceci permet de d\'efinir {\it l'\'etat spinoriel} $sp(C)$ de chaque courbe rationnelle r\'eelle $C$ comme valant $+1$ dans
le premier cas, et $-1$ dans le second. Lorsque les lieux r\'eels de $L$ et $M$ ne sont pas orientables, on modifie 
ces axes mobiles \`a l'aide d'un demi-tour \`a droite donn\'e par l'orientation de $\R X$, ce qui permet de se ramener au cas pr\'ec\'edent et de
d\'efinir l'\'etat spinoriel \'egalement dans ce cas.
L'entier $\chi^{d, \s}_r (\underline{x})$ n'est alors autre que le nombre de courbes rationnelles 
r\'eelles qui r\'ealisent la classe d'homologie $d$ et passent par $\underline{x}$, ces courbes \'etant compt\'ees en fonction
de leur \'etat spinoriel, de sorte que
$$\chi^{d,\s}_{{r}} (\underline{x}) = \sum_{C \in {\cal R}_d (\underline{x})} sp(C) \in \Z.$$
On a not\'e $r = (r_1 , \dots , r_N)$ le $N$-uplet associ\'e \`a $\underline{x}$ ; c'est-\`a-dire que $N$ d\'esigne le nombre de composantes connexes de $\R X$
et en notant  $(\R X)_1 , \dots , (\R X)_N$ ces composantes, $r_i = \# (\underline{x} \cap (\R X)_i)$.

\begin{theorem}[\cite{WelsDuke}]
\label{theoDuke}
Soient $(X , c_X)$ une vari\'et\'e alg\'ebrique r\'eelle convexe lisse de dimension trois et $\s$ une structure $Spin_3$ 
sur son lieu r\'eel $\R X$ suppos\'e orientable. Soit $d \in H_2 (X ; \Z)$ telle que $c_1 (X)d$ soit pair et diff\'erent de quatre et soit
$k_d = \frac{1}{2} c_1 (X)d \in \N^*$. Soient $\underline{x}= (x_1 , \dots , x_{k_d})$ une configuration r\'eelle
de $k_d$ points distinct dont au moins un r\'eel et $r$ le $N$-uplet associ\'e.  
L'entier $\chi^{d,\s}_{\underline{r}} (\underline{x})$ ne d\'epend alors pas du choix g\'en\'erique de $\underline{x}$.
\end{theorem}

Ce r\'esultat est valable aussi pour les vari\'et\'es dont le lieu r\'eel n'est pas orientable, moyennant le choix d'une structure $Pin_3^-$ sur le lieu r\'eel ; l'invariant
s'annule alors lorsque $k_d$ est pair, voir \cite{WelsDuke}.

Le Th\'eor\`eme \ref{theoDuke} permet sans ambigu\"{\i}t\'e de noter cet entier $\chi^{d, \s}_r$. Lorsque $\sum_{i=1}^N r_i$ n'a pas la m\^eme parit\'e que $k_d = \frac{1}{2} c_1 (X) d$, 
on pose $\chi^{d, \s}_r = 0$.
On note alors $\chi^{d, \s} [T]$ la fonction g\'en\'eratrice $\sum_{|r|=0}^{k_d} \chi_r^{d, \s} T^r \in \Z [T_1 , \dots , T_N]$.
Cette fonction est de m\^eme parit\'e que $\frac{1}{2} c_1 (X) d$ et tous ses mon\^omes ne d\'ependent en fait que d'une ind\'etermin\'ee. 

Ainsi, la fonction $\chi^{\s} : d \in H_2 (X ; \Z) \mapsto \chi^{d,\s} (T) \in \Z[T]$ est invariante par isomorphisme de la vari\'et\'e alg\'ebrique r\'eelle convexe lisse de dimension trois $(X , c_X)$.
On en d\'eduit \`a nouveau les bornes inf\'erieures suivantes en g\'eom\'etrie \'enum\'erative r\'eelle.

\begin{cor}[\cite{WelsDuke}]
\label{corlowerDuke}
Sous les hypoth\`eses du Th\'eor\`eme \ref{theoDuke}, notons ${R}_d (\underline{x})$
le nombre de courbe rationnelles r\'eelles connexes homologues \`a $d$ qui passent par $\underline{x}$ 
et $N_d$ l'invariant de Gromov-Witten de genre z\'ero associ\'e. Alors,
$|\chi^{d,\s}_r| \leq R_d (\underline{x}) \leq N_d. \quad \square$
\end{cor}

Finissons ce paragraphe par une interpr\'etation topologique de nos r\'esultats.
Les singularit\'es de l'espace $\R_\tau \overline{\cal M}^d_{k_d} (X)$ sont de codimension au moins deux, de sorte que cet espace
poss\`ede une premi\`ere classe de Stiefel-Whitney. \'Etant donn\'e $D \in H_{3k_d - 1} (\R_\tau \overline{\cal M}^d_{k_d} (X) ; \Z/2\Z)$, 
on note $D^\vee$ son image sous le morphisme $H_{3k_d - 1} (\R_\tau \overline{\cal M}^d_{k_d} (X) ; \Z/2\Z)
\to H^1 (\R_\tau \overline{\cal M}^d_{k_d} (X) ; \Z/2\Z)$.
\begin{prop}[\cite{WelsSemi}]
\label{propw1}
La  premi\`ere classe de Stiefel-Whitney de toute composante $\R {\cal M}^*$ de $\R_\tau {\cal M}^d_{k_d} (X)$ 
qui contient une courbe \'equilibr\'ee s'\'ecrit 
$$w_1 (\R {\cal M}^*) = (\R_\tau ev^d)^* w_1 (\R_\tau X^{k_d}) + 
\sum_{D \subset \text{Red}'} \epsilon (D) D^\vee \in H^1 (\R {\cal M}^* ; \Z/2\Z),$$
o\`u $\epsilon (D) \in \{ 0, 1 \}$ et lorsque $\epsilon (D) = 1$, la composante irr\'eductible $D$ de
$\text{Red}$ se trouve contract\'ee par l'application d'\'evaluation $\R_\tau ev^d$. $\square$
\end{prop}
On a not\'e ici $\text{Red}'$ la r\'eunion du diviseur des courbes r\'eductibles $\text{Red}$ et de l'\'eventuel diviseur des courbes
non-\'equilibr\'ees $(u , C , \underline{z})$ telles que $\dim H^1 (C ; N_u \otimes
{\cal O}_C (-  \underline{z})) \geq 2$, si un tel diviseur existe. On note $\text{Red}_1$ 
la r\'eunion des composantes irr\'eductibles $D$ de $\text{Red}'$ pour lesquelles $\epsilon (D) = 1$.
En dimension deux, cet ensemble a \'et\'e d\'etermin\'e dans \cite{Pui}. 
\'Equipons $\R_\tau X^{k_d}$ d'un syst\`eme de coefficients tordus entiers ${\cal Z}$  et notons 
$[\R_\tau X^{k_d}] \in H_{3k_d} (\R_\tau X^{k_d} ; {\cal Z})$ sa classe fondamentale. Notons 
${\cal Z}^*$ le syst\`eme de coefficients locaux induit sur $\R {\cal M}^*$, tir\'e en arri\`ere de ${\cal Z}$ par $\R_\tau ev^d$. 
\begin{prop}[\cite{WelsSemi}]
\label{propfond}
Sous les hypoth\`eses de la Proposition \ref{propw1}, il existe une unique classe fondamentale 
$[\R {\cal M}^*] \in H_{3k_d} (\R {\cal M}^* , \text{Red}_1 ; {\cal Z}^*)$ telle qu'en toute courbe \'equilibr\'ee
$(u, C , \underline{z}) \in \R {\cal M}^*$, le morphisme
$$(\R_\tau ev^d)_* : H_{3k_d} (\R {\cal M}^* , 
\R {\cal M}^* \setminus \{ (u, C , \underline{z}) \} ; {\cal Z}^*) \to
H_{3k_d} (\R_\tau X^{k_d} , \R_\tau X^{k_d} \setminus \{ u(\underline{z}) \} ; {\cal Z})$$
envoie
$[\R {\cal M}^*]$ sur $sp (u, C , \underline{z}) [\R_\tau X^{k_d}]$. $\square$
\end{prop}
Comme $\R_\tau ev^d (\text{Red}_1)$ est de codimension deux, le groupe $H_{3k_d} (\R_\tau X^{k_d} , 
\R_\tau ev^d (\text{Red}_1) ; {\cal Z})$ est cyclique, engendr\'e par $[\R_\tau X^{k_d}]$. 
L'entier 
$\chi_r^{d, \s}$ n'est autre que celui d\'efini par la relation $(\R_\tau ev^d)_* 
[\R_\tau \overline{\cal M}^d_{k_d} (X)] = \chi_r^{d, \s} [\R_\tau X^{k_d}]$, o\`u la classe fondamentale $[\R_\tau \overline{\cal M}^d_{k_d} (X)]$ 
est donn\'ee par la Proposition \ref{propfond}.

\subsection{Extension aux vari\'et\'es symplectiques r\'eelles fortement semi-positives}

L'extension des r\'esultats du \S \ref{subsectDuke} aux vari\'{e}t\'{e}s symplectiques n'est pas imm\'ediate, en partie parce-que le th\'eror\`eme de Grothendieck \cite{Grot} selon lequel les fibr\'es holomorphes
sur la sph\`ere de Riemann sont enti\`erement d\'ecomposables n'est plus valable pour les fibr\'es normaux des courbes pseudo-holomorphes. Ces derniers sont des
fibr\'es vectoriels complexes munis d'un op\'erateur de Cauchy-Riemann qui n'est que $\R$-lin\'eaire et non $\C$-lin\'eaire comme dans le cas de fibr\'es
holomophes. Ces premiers sont des perturbations d'ordre z\'ero de ces derniers par des op\'erateurs $\C$-antilin\'eaires et sont parfois appel\'es 
\og op\'erateurs de Cauchy-Riemann g\'en\'eralis\'es \fg.  J'ai \'etendu dans  \cite{WelsSemi} la notion d'\'etat spinoriel pour un op\'erateur de Cauchy-Riemann g\'en\'eralis\'e surjectif.

La strat\'egie est la suivante. L'espace des op\'erateurs de Cauchy-Riemann g\'en\'eralis\'es r\'eels sur un fibr\'e vectoriel complexe r\'eel donn\'e est un espace de Banach affine, il contient
les op\'erateurs de Cauchy-Riemann $\C$-lin\'eaires comme sous-espace de Banach. Or chaque op\'erateur de Cauchy-Riemann surjectif d\'efinit une structure de fibr\'e vectoriel holomorphe
\'equilibr\'e et poss\`ede donc un \'etat spinoriel d'apr\`es les r\'esultats du \S \ref{subsectDuke}. \'Etant donn\'e un op\'erateur de Cauchy-Riemann g\'en\'eralis\'e surjectif, on le relie \`a
op\'erateur de Cauchy-Riemann surjectif par un chemin g\'en\'erique et on d\'efinit son \'etat spinoriel comme celui de l'op\'erateur de Cauchy-Riemann si le chemin traverse un nombre pair
de fois le mur des op\'erateurs non-surjectifs et son oppos\'e sinon. 

Soit alors $(X , \omega , c_X)$ une vari\'et\'e symplectique r\'eelle fortement semi-positive de dimension six, c'est-\`a-dire pour laquelle toute classe sph\'erique
$d \in H_2 (X ; \Z)$ positive contre $\omega$ satisfait l'implication $c_1(X) d \geq 2 - n \implies c_1(X) d \geq 1$. Les vari\'et\'es symplectiques r\'eelle positives, 
par exemple de Fano, satisfont cette condition. On suppose \`a nouveau pour simplifier le lieu r\'eel de cette vari\'et\'e orientable et on l'\'equipe d'une structure spin $\s$.
Soit,  comme au  \S \ref{subsectDuke}, 
$d \in H_2 (X ; \Z)$ telle que $(c_X)_* d = -d$, $c_1(X)d $ soit pair et strictement plus grand que deux. 
Soient $k_d = \frac{1}{2} c_1(X)d $ et  $\underline{x} = (x_1 , \dots , 
x_{k_d}) \in X^{k_d}$ une configuration r\'eelle de $k_d$ points distincts, dont au moins un r\'eel. 
Lorsque $J \in \R {\cal J}_\omega$ est suffisamment g\'en\'erique, il n'y a qu'un nombre fini de courbes $J$-holomorphes rationnelles r\'eelles connexes
homologues \`a $d$ et contenant $\underline{x}$. Ces courbes sont toutes irr\'eductibles, lisses et de partie r\'eelle non-vide. On
note ${\cal R}_d (\underline{x} , J )$ cet ensemble fini de courbes.
Le fibr\'e normal de chacune de ces courbes $C \in {\cal R}_d (\underline{x} , J )$
est \'equip\'e d'un op\'erateur de Cauchy-Riemann g\'en\'eralis\'e surjectif $D_C$ qui poss\`ede donc un \'etat spinoriel $sp (C)$ d'apr\`es ce qui pr\'ec\`ede.
On pose
$$\chi_r^{d,  \s} (\underline{x} , J) = \sum_{C \in {\cal R}_d (\underline{x} , J)}
sp (C) \in \Z.$$
\begin{theorem}[\cite{WelsSemi}]
\label{theoSemi}
Soit $(X , \omega , c_X)$ une vari\'et\'e symplectique r\'eelle fortement semi-positive  de dimension six, de lieu r\'eel orientable muni d'une structure spin $\s$.
Soit $d \in H_2 (X ; \Z)$ telle que $(c_X)_* d = -d$, $c_1(X)d $ est pair et strictement plus grand que deux. 
Soient $k_d = \frac{1}{2} c_1(X)d$ et  $\underline{x}$ une configuration r\'eelle de $k_d$ points distincts, dont au moins un r\'eel. 
et $r = (r_1 , \dots , r_N)$ le $N$-uplet associ\'e. 
Alors, l'entier $\chi_r^{d, \s} (\underline{x} , J)$ ne d\'epend ni du choix de  $\underline{x}$, ni du choix g\'en\'erique de $J \in \R {\cal J}_\omega$.
\end{theorem}

Remarquons que ce r\'esultat 
permet de noter sans ambigu\"{\i}t\'e l'invariant $\chi_r^{d, \s}$, c'est un invariant par d\'eformation fortement semipositive de $(X , \omega , c_X)$.
On en d\'eduit les bornes inf\'erieures suivantes.
\begin{cor}[\cite{WelsSemi}]
\label{corlowerboundsSemi}
Sous les hypoth\`eses du Th\'eor\`eme \ref{theoSemi},
$|\chi_r^{d,  \s}| \leq \# {\cal R}_d (\underline{x} , J),$
pour tout choix de configuration r\'eelle $\underline{x} \in X^{k_d}$ telle que $\underline{x} \cap \R X = r$,
et tout choix g\'en\'erique de $J \in \R {\cal J}_\omega$. $\square$
\end{cor}

Les invariants qui ressortent des Th\'eor\`emes \ref{theoInvent} et \ref{theoDuke} ont \'et\'e interpr\'et\'es par C.-H. Cho \cite{Cho}
et J. Solomon \cite{Sol}. Leur approche consiste \`a d'abord d\'efinir la classe fondamentale $[\R_\tau \overline{\cal M}^d_{k_d} (X)]$ 
donn\'ee par la Proposition \ref{propfond} en utilisant les travaux de K. Fukaya, Y.-G. Oh, H. Ohta et K. Ono \cite{FOOOI},  \cite{FOOOII},
puis \`a en d\'eduire l'existence des invariants gr\^ace \`a la relation entre classes fondamentales donn\'ee \`a la suite de cette proposition. 
J. Solomon a \'etendu ces invariants aux courbes de genre strictement positifs
mais de structure conforme fix\'ee et aux vari\'et\'es symplectiques de dimension six, notamment de Calabi-Yau. Dans le cas des quintiques de $\C P^4$, l'invariant a \'et\'e
 calcul\'e par R. Pandharipande, J. Solomon et J. Walcher \cite{PSW}.

\subsection{Optimalit\'e, congruences et calculs dans le cas de l'ellipso\"{\i}de de dimension trois}

\begin{theorem}[\cite{WelsSFT}]
\label{theoopt4} 
Soient $(X, c_X)$ la quadrique ellipso\"{\i}de de dimension trois et  $d \in H_2 (X ; \Z)$ satisfaisant $c_1 (X) d = 2 \mod (4)$. L'invariant $\chi_1^d$ est alors n\'egatif et les bornes inf\'erieures 
apparues dans le Corollaire \ref{corlowerDuke} sont optimales, atteintes lorsque les conditions d'incidence non r\'eelles sont choisies suffisamment proches d'une section 
 hyperplane r\'eelle disjointe du lieu r\'eel $\R X$.
\end{theorem}

\begin{rem}
La condition $c_1 (X) d = 2 \mod (4)$ garantit la parit\'e de l'entier $k_d$ de sorte que l'on peut effectivement choisir un point r\'eel. Lorsque $c_1 (X) d = 0 \mod (4)$,
et $r = 0$, l'invariant $\chi_r^d$ n'est pas d\'efini. Toutefois, on a montr\'e dans ce cas l\`a qu'il existe une structure presque-complexe g\'en\'erique $J \in \R {\cal J}_\omega$
et $k_d$ points complexes conjugu\'es pour lesquels aucune courbe $J$-holomorphe rationnelle r\'eelle homologue \`a $d$ contient ces $k_d$ points, voir le Th\'eor\`eme \ref{theomin2}.
\end{rem}

\begin{theorem}[\cite{WelsSFT}]
\label{theocong2}
Soient $(X, c_X)$ la quadrique ellipso\"{\i}de de dimension trois et $d$ un multiple positif, disons $\delta > 0$, d'une section hyperplane r\'eelle.
Lorsque $6r + 1 \leq 3 \delta$, la puissance $2^{\frac{3}{4} (\delta - 2r)}$ divise $\chi^d_r$. 
\end{theorem}

\begin{cor}[\cite{WelsSFT}]
\label{corcalc3spher}
Soit $(X, \omega , c_X)$ une vari\'et\'e symplectomorphe \`a la quadrique ellipso\"{\i}de de dimension trois. Alors,
$\chi^{2}_1 = - 1$, $\chi^{6}_1 = 0$ et $\chi^{10}_1 = -896$.
\end{cor}

Dans le cas de l'espace projectif de dimension trois, une formule calculant $\chi^d_{2d} (\C P^3)$ pour tout degr\'e $d$ est annonc\'ee par
E. Brugall\'e et G. Mikhalkin dans \cite{BruMikh}. En particulier,  $\chi^5_{10} = 45$ et $\chi^7_{14} = -14589$ alors qu'en degr\'e pair, l'invariant s'annule pour
des raisons de sym\'etrie.

\section{Sur la pr\'esence et l'absence de membranes $J$-holomorphes}

\subsection{Absence de membranes $J$-holomorphes}

Soit $C$ une membrane $J$-holomorphe \`a bord dans une sous-vari\'et\'e lagrangienne $L$ d'une vari\'et\'e symplectique ferm\'ee $(X, \omega)$. Notons $\chi$ la caract\'eristique d'Euler de
cette membrane, $d \in H_2 (X , L ; \Z)$ sa classe d'homologie relative et $\mu_{TX} \in H^2 (X , L ; \Z)$ la classe de Maslov de la paire $(X , L)$. La dimension attendue de l'espace des d\'eformations
de $C$ s'\'ecrit $\langle \mu_{TX} , d \rangle + (n-3)\chi$. Cette dimension chute lorsque l'on impose \`a $C$ des contraintes suppl\'ementaires. Si l'on impose par exemple
\`a cette membrane de rencontrer $p$ cycles de codimensions $2 + q_1, \dots , 2 + q_p$, cette dimension attendue chute de la somme $q = q_1+ \dots + q_p$. Deux probl\`emes g\'en\'eraux sous-tendent
nos r\'esultats. Il s'agit d'une part de compter les membranes $J$-holomorphes homologues \`a $d$ soumises \`a de telles conditions d'incidence de sorte que ce comptage ne d\'epende pas
de $J$ et ne d\'epende des conditions d'incidence qu'\`a homologie pr\`es. Il s'agit d'autre part de minimiser ce nombre de membranes. Si nous ne pouvons r\'epondre au premier probl\`eme
dans ce degr\'e de g\'en\'eralit\'e, il nous est par contre parfois possible de r\'epondre au second sans m\^eme supposer l'\'egalit\'e $q = \langle \mu_{TX} , d \rangle + (n-3)\chi$, lorsque le minimum en question est nul. Le pr\'esent paragraphe est consacr\'e aux  r\'esultats que l'on a pu obtenir
dans cette direction. Ici encore le minimum est atteint en allongeant le cou d'une structure presque complexe g\'en\'erale.

\subsubsection{En dimension sup\'erieure}

\begin{theorem}[\cite{WelsSFT}, \cite{WelsFloer}]
\label{theomin3} 
Soit $L$ une sph\`ere lagrangienne dans une vari\'et\'e symplectique ferm\'ee $(X, \omega)$ satisfaisant $c_1 (X) = \lambda \omega$, $\lambda \leq 0$
et soit $E > 0$. Supposons la dimension de $X$ sup\'erieure \`a cinq. Pour toute structure presque-complexe $J$ g\'en\'erale ayant un cou suffisamment long au voisinage de $L$,
cette vari\'et\'e ne poss\`ede ni membrane $J$-holomorphe reposant sur $L$ ni courbe $J$-holomorphe rencontrant $L$ qui soit d'\'energie inf\'erieure \`a $E$. 
Ce r\'esultat reste valable en dimension quatre pour les courbes ou membranes de genre nul. 
\end{theorem}
Rappelons que l'\'energie d'une courbe $C$ est par d\'efinition l'int\'egrale de la forme $\omega$ sur cette courbe. Les vari\'et\'es projectives \`a fibr\'e canonique nul ou ample, par 
exemple les intersections compl\`etes de multidegr\'es $(d_1, \dots , d_k)$
de l'espace projectif de dimension $N$ d\`es lors que $\sum_{i=1}^k d_i \geq N+1$, satisfont les hypoth\`eses du Th\'eor\`eme \ref{theomin3}. Remarquons qu'une modification de ce dernier
s'applique \'egalement aux vari\'et\'es dont le fibr\'e canonique est le produit d'un fibr\'e ample et d'un fibr\'e port\'e par un diviseur effectif disjoint de $L$. Le Th\'eor\`eme \ref{theomin3} permet de d\'efinir la cohomologie de Floer de sph\`eres lagrangiennes dans les vari\'et\'es symplectiques dont la premi\`ere classe de Chern s'annule, voir \cite{WelsFloer}
et \cite{FOOOI}, \cite{FOOOII} pour une th\'eorie de l'obstruction \`a d\'efinir en g\'en\'eral une telle homologie.

\begin{theorem}[\cite{WelsSFT}]
\label{theomin2} 
Soit $L$ une sph\`ere lagrangienne dans une vari\'et\'e symplectique ferm\'ee semipositive $(X, \omega)$ de dimension $2n \geq 6$ et soit $d \in H_2 (X , L ; \Z)$. \'Ecrivons 
$\langle \mu_{TX} , d \rangle + (n-3)\chi = q + r$ avec $q \in \Z$, $0 \leq r < 2 + (n-3)\chi$ et $\chi \leq 2$. Lorsque $q \geq 0$, choisissons $p$ cycles de $X \setminus L$ de codimensions 
$2 + q_1, \dots , 2 + q_p$ de sorte que $q = q_1+ \dots + q_p$.  D\`es que la structure presque complexe g\'en\'erale $J$ poss\`ede un cou suffisamment long au voisinage 
 de $L$, cette vari\'et\'e ne contient aucune membrane $J$-holomorphe homologue \`a $d$, de caract\'eristique d'Euler $\chi$ qui rencontre ces $p$ cycles et repose sur $L$.
Ce r\'esultat reste valable pour des membranes de genre nul lorsque $n = 2$.
\end{theorem}

{\bf Exemple : la quadrique ellipso\"{\i}de.}

Soit $X$ la quadrique ellipso\"{\i}de de dimension complexe $n \geq 3$ et $H$ une section hyperplane disjointe de $L$. Le groupe $H_2 (X , L ; \Z)$ est monog\`ene, engendr\'e par
la classe $d_0$ satisfaisant $\langle H , d_0 \rangle = +1$. La premi\`ere classe de Chern de $X$ vaut $n H$, d'o\`u l'on d\'eduit le calcul $\langle \mu_{TX} , l d_0 \rangle
= 2ln$ quel que soit l'entier $l$. \'Ecrivons $l = (n-1)a + b$, le Th\'eor\`eme \ref{theomin2} s'applique par exemple lorsque $n+1 \leq 2b < 2n$, les membranes sont des disques
et lorsque toutes les conditions d'incidence sont ponctuelles.
Rappelons que le Th\'eor\`eme \ref{theoopt4} traite du cas $r=n-1$ et montre ainsi en un sens l'optimalit\'e des hypoth\`eses faites dans ce Th\'eor\`eme \ref{theomin2}.

\subsubsection{En dimension quatre}

Nous noterons ${\cal M}_{g,b}$ l'espace des modules des structures complexes de la surface compacte connexe orient\'ee de genre $g$ ayant $b$ composantes de bord.

\begin{prop}[\cite{WelsSFT}]
\label{propminsphere}
Soit $L$ une sph\`ere lagrangienne dans une vari\'et\'e symplectique ferm\'ee de dimension quatre $(X, \omega)$. On suppose que cette derni\`ere ne poss\`ede pas de sph\`ere
symplectique $S$ satisfaisant $\langle c_1 (X) , [S] \rangle > 0$. Soit $(d, g, b) \in H_2 (X , L ; \Z) \times \N \times \N^*$ et $K$ un compact de ${\cal M}_{g,b}$. Alors, pour toute
structure presque-complexe g\'en\'erale ayant un cou suffisamment long au voisinage de $L$, la vari\'et\'e ne poss\`ede pas de membrane $J$-holomorphe homologue \`a $d$ \`a 
bord dans $L$ et conforme \`a un \'el\'ement de $K$.
\end{prop}

\begin{prop}[\cite{WelsSFT}]
\label{propmin1} 
Soit $L$ une surface lagrangienne orientable hyperbolique dans une vari\'et\'e symplectique ferm\'ee de dimension quatre $(X, \omega)$ et  soit $d \in H_2 (X , L ; \Z)$. On note
$N_d^g (\underline{x} , J)$ le nombre de courbes $J$-holomorphes homologues \`a $d$ \`a bords dans $L$, de topologie et de structure conforme donn\'ees et qui passent par
une configuration $\underline{x}$ de points distincts de $(X, \omega)$ de cardinal ad\'equat, pour $J \in {\cal J}_\omega$ g\'en\'erique. Ce nombre
$N_d^g (\underline{x} , J)$ s'annule pour toute structure presque-complexe g\'en\'erale ayant un cou suffisamment long au voisinage de $L$.
\end{prop}

\begin{prop}[\cite{WelsSFT}]
\label{propmintore}
Soit $(X, \omega , c_X)$ une vari\'et\'e symplectique r\'eelle ferm\'ee de dimension quatre dont le lieu r\'eel poss\`ede un tore lagrangien ou bien une surface hyperbolique lagrangienne $L$,
orientable ou non. On suppose que  $(X, \omega , c_X)$ ne poss\`ede pas de sph\`ere symplectique r\'eelle $S$ satisfaisant $\langle c_1 (X) , [S] \rangle > 1$ si $L$ est orientable et 
$\langle c_1 (X) , [S] \rangle > 0$ sinon. Soit $(d, g, b) \in H_2 (X , L ; \Z) \times \N \times \N^*$ 
et $K$ un compact de ${\cal M}_{g,b}$. Alors, pour toute structure presque-complexe g\'en\'erale ayant un cou suffisamment long au voisinage de $L$, la vari\'et\'e ne poss\`ede pas de membrane 
$J$-holomorphe homologue \`a $d$ \`a bord dans $L$ et conforme \`a un \'el\'ement de $K$.
\end{prop}

\subsection{Pr\'esence de membranes $J$-holomorphes}
\label{subsectpre}

Les r\'esultats pr\'esent\'es aux \S \S \ref{subsectinv} et \ref{subsectDuke} permettent de garantir l'existence de disques $J$-holomorphes reposant sur une sous-vari\'et\'e lagrangienne d'une
vari\'et\'e symplectique donn\'ee, lorsque cette lagrangienne se trouve dans le lieu fixe d'une involution antisymplectique, laquelle est $J$-antiholomorphe et \`a condition que l'invariant
que l'on a d\'efini n'est pas nul. Nous souhaitons montrer ici qu'il est possible d'obtenir ces r\'esultats pour une classe plus large de sous-vari\'et\'e lagrangiennes, en faisant intervenir
la notion d'involutions antibirationnelles sur les vari\'et\'es symplectiques. 

\subsubsection{Involutions antibirationnelles des vari\'et\'es symplectiques de dimension quatre}

Une involution $c_X$ de la vari\'et\'e symplectique de dimension quatre $(X , \omega)$ qui est d\'efinie en-dehors d'un nombre fini de points
$x_1 , \dots , x_k$ de $X$ est dite {\it antibirationnelle} lorsqu'il existe un diagramme commutatif de la forme suivante :
$$\begin{array}{rcl}
\label{diagram}
(Y, J_Y) & \stackrel{c_Y}{\longrightarrow} & (Y, J_Y)\\
\pi \downarrow && \downarrow \pi \\
(X, J_X) & \stackrel{c_X}{\longrightarrow} & (X, J_X)\\
\end{array}$$
o\`u $Y$ est une vari\'et\'e compacte de dimension quatre obtenue \`a partir de $X$ en r\'ealisant un
nombre fini d'\'eclatements topologiques au-dessus des points $x_i$, $i \in \{1, \dots , k\}$, $J_X , J_Y$ sont
des structures presque-complexes lisses et $c_Y$ une involution $J_Y$-antiholomorphe sur $Y$ toute enti\`ere. 
De plus, $J_X$ est suppos\'ee $\omega$-positive, $c_X$ est
 $J_X$-antiholomorphe sur son lieu de d\'efinition et  $\pi$ est $(J_Y , J_X)$-holomorphe.

Les involutions  antibirationnelles classiques sur les surfaces compactes de K\"ahler fournissent des exemples de telles surfaces.
Remarquons que pour tout $i \in \{1, \dots , k\}$,
$\pi^{-1} (x_i)$ est un arbre de sph\`eres $J_Y$-holomorphes n'ayant que des points doubles transverses comme singularit\'es. 

\begin{lemma}
\label{lemmaassume}
Supposons que pour tout $i \in \{1, \dots , k\}$ et toute composante irr\'eductible $C$ de l'arbre
$\pi^{-1} (x_i)$, $c_Y (C)$ ne soit pas contract\'ee par $\pi$ sur $x_1 , \dots , x_k$.
Alors, le diagramme ci-dessus est unique \`a \'equivalence pr\`es, une fois donn\'ee $(X, \omega , c_X)$.
\end{lemma}

Soient $(X, \omega , J_X , c_X)$ satisfaisant les hypoth\`eses du Lemme \ref{lemmaassume} et
$(Y, J_Y ,c_Y)$ la vari\'et\'e de dimension quatre associ\'ee. Soit
$\underline{y}$ l'ensemble fini $\big( \cup_{i=1}^k \pi^{-1} (x_i) \big) 
\cap c_{Y} \big( \cup_{i=1}^k \pi^{-1} (x_i) \big)$. L'involution antibirationnelle
$c_X$ est dite {\it simple} lorsqu'elle satisfait les hypoth\`eses du Lemme \ref{lemmaassume} et lorsque
$\underline{y}$ se trouve en-dehors des points doubles de $\cup_{i=1}^k \pi^{-1} (x_i)$.

\begin{lemma}
\label{lemma2form}
Soit $c_X$ une involution antibirationnelle simple de  $(X , \omega)$ et
$(Y, c_Y)$ la vari\'et\'e de dimension quatre donn\'ee par le Lemme \ref{lemmaassume}. Alors,
la deux-forme $\omega_Y = \pi^* \omega - (\pi \circ c_Y)^* \omega$ est ferm\'ee et non-d\'eg\'en\'er\'ee 
en tout point de $Y \setminus \underline{y}$. Elle est \'egalement non-d\'eg\'en\'er\'ee 
en tout point d'intersection transverse de $\big( \cup_{i=1}^k 
\pi^{-1} (x_i) \big) \cap c_{Y} \big( \cup_{i=1}^k \pi^{-1} (x_i) \big)
\subset \underline{y}$.
\end{lemma}

Une telle deux-forme qui n'a qu'un nombre fini de noyaux de dimension deux sera dite {\it quasi-symplectique}.
Remarquons qu'en particulier, lorsque l'intersection $\big( \cup_{i=1}^k 
\pi^{-1} (x_i) \big) \cap c_{Y} \big( \cup_{i=1}^k \pi^{-1} (x_i) \big)$ est
transverse, la deux-forme $\omega_Y$ est symplectique.\\

La structure presque-complexe $J_Y$ est $\omega_Y$-positive dans le sens que pour tous $y \in Y$ et $v \in T_y Y \setminus \{ 0 \}$, soit $v$ et
$J_Y (v)$ engendrent le noyau de $\omega_Y|_y$, soit $\omega_Y (v , J_Y (v))>0$. Notons ${\cal J}_{\omega_Y}$ l'espace des structures presque-complexes de classe $C^l$
qui sont $\omega_Y$-positives.  Si $J \in {\cal J}_{\omega_Y}$,
alors $\overline{c}_Y^* (J) = -dc_Y \circ J \circ dc_Y$ appartient \'egalement \`a
${\cal J}_{\omega_Y}$. Notons $\R {\cal J}_{\omega_Y}$ le lieu fixe de cette action de
$\Z/2 \Z$ sur ${\cal J}_{\omega_Y}$.

\begin{lemma}
\label{lemmaJ0}
Soient $c_X$ une involution antibirationnelle simple sur $(X , \omega)$ et
$(Y, c_Y ,\omega_Y )$ la vari\'et\'e de dimension quatre donn\'ee par les Lemmes \ref{lemmaassume},
\ref{lemma2form}. Il existe $J_0 \in \R {\cal J}_{\omega_Y}$ tel que $\omega_Y (J_0 , J_0) =
\omega_Y$ et $g_Y = \omega_Y (. , J_0 )$ soit un deux-tenseur sym\'etrique 
positif sur $Y$, d\'efini en-dehors de $\underline{y}$.
\end{lemma}

Pour tout voisinage $U$ de $\underline{y}$ et tout $J_0 \in \R {\cal J}_{\omega_Y}$
donn\'e par le Lemme \ref{lemmaJ0}, notons ${\cal J}_{\omega_Y}^{U , J_0}$ (resp.
$\R {\cal J}_{\omega_Y}^{U , J_0}$) le sous-espace des $J \in 
{\cal J}_{\omega_Y}$ (resp. $J \in \R {\cal J}_{\omega_Y}$)
telles que $J = J_0$ sur $U$.

\begin{lemma}
\label{lemmaJ0}
Pour tous $U , J_0$, l'espace ${\cal J}_{\omega_Y}^{U , J_0}$ est une vari\'et\'e de Banach s\'eparable non-vide
et contractile. Le sous-espace $\R {\cal J}_{\omega_Y}^{U , J_0}$ en est une sous-vari\'et\'e de Banach s\'eparable non-vide
et contractiles.$\square$ 
\end{lemma}

\begin{rem}
La deux-forme
$\pi^* \omega$ est limite d'une suite de formes symplectiques sur $Y$ obtenues
apr\`es un nombre fini d'\'eclatements de boules symplectiques dont les rayons convergent vers z\'ero.
Par suite, la deux-forme $\omega_Y$ est limite d'une suite de formes symplectiques 
$(\omega_Y^n)_{n \in \N}$. Alors,  $J \in {\cal J}_{\omega_Y}$ est $\omega_Y^n$-positif pour $n$ assez grand, principalement parce-que les noyaux de $\omega_Y$ 
deviennent 
symplectiques pour $\omega_Y^n$.
\end{rem}

\subsubsection{Invariants \'enum\'eratifs des involutions antibirationnelles simples}

Soient $c_X$ une involution antibirationnelle simple sur $(X , \omega)$ et
$x_1 , \dots , x_k \in X$ les points o\`u elle n'est pas d\'efinie. Soit
$(Y , \omega_Y , c_Y)$ la vari\'et\'e quasi-symplectique de dimension quatre associ\'ee, voir le
Lemme \ref{lemma2form}. Soient $\pi$ la projection $Y \to X$
et $\underline{y}$ l'ensemble fini $\big( \cup_{i=1}^k \pi^{-1} (x_i) \big) 
\cap c_{Y} \big( \cup_{i=1}^k \pi^{-1} (x_i) \big)$. Soit $\R Y$ le lieu fixe de
$c_Y$, on \'etiquette ses composantes connexes $(\R Y)_1, \dots , (\R Y)_N$. Remarquons
que la courbe $\big( \cup_{i=1}^k \pi^{-1} (x_i) \big) 
\cup c_{Y} \big( \cup_{i=1}^k \pi^{-1} (x_i) \big)$ n'intersecte $\R Y$ qu'en un nombre fini de
points, de sorte qu'elle ne d\'{e}connecte aucune des courbes $(\R Y)_i$, 
$i \in \{1, \dots , N\}$. Soient $d_Y \in H_2 (Y  ; \Z)$ tel que $(c_Y)_* d_Y = -d_Y$,
$c_1 (Y)d_Y > 0$ et $y = (y_1 , \dots y_{c_1 (Y)d_Y - 1})$ une configuration
r\'eelle de $c_1 (Y)d_Y - 1$ points distinct de $Y \setminus \big( \cup_{i=1}^k \pi^{-1} (x_i) 
\cup c_{Y} \big( \cup_{i=1}^k \pi^{-1} (x_i) \big) \big)$. Pour tout $i \in \{1 , \dots , N\}$,
notons $r_i = \# (y \cap (\R Y)_i)$ puis $r = (r_1 , \dots , r_N)$. Soient
$U$, voisinage de $\underline{y}$ et $J_0 \in \R {\cal J}_{\omega_Y}$ donn\'es par le
 Lemme \ref{lemmaJ0}. Alors, d\`es que $U$ est suffisamment petit, 
pour tout $J \in \R {\cal J}_{\omega_Y}^{U , J_0}$ g\'en\'erique, il n'y a qu'un nombre fini
 de courbes $J$-holomorphes rationnelles r\'eelles homologue \`a $d_Y$ dans $Y$ qui contiennent $y$. 
 Ces courbes sont toutes irr\'eductibles, immerg\'ees et n'ont que des points doubles transverses comme singularit\'es.
 Le nombre total de leurs points doubles vaut $\delta_Y = \frac{1}{2} (d_Y^2 - c_1 (Y) d_Y +2)$. 
Pour tout entier $m$ compris entre $0$ et $\delta_Y$, notons $n_d (m)$ le nombre de ces courbe qui sont de masse $m$. On pose alors
$$\chi_r^{d_Y} (y,J, U , J_0) = \sum_{m=0}^{\delta_Y} (-1)^m n_d (m).$$
\begin{theorem}
\label{theoantibir}
L'entier $\chi_r^{d_Y} (y,J , U , J_0)$ ne d\'epend pas des choix de $y,J , U$ et $J_0$. $\square$
\end{theorem}

Remarquons que l'entier 
$\chi_r^{d_Y}$ fourni par le Th\'eor\`eme \ref{theoantibir} est un invariant par d\'eformation du triplet 
 $(X , \omega , c_X)$, puisque le triplet $(Y , \omega_Y , c_Y)$ lui est canoniquement associ\'e. 

\subsubsection{Exemple: les tores  isotopes au tore de Clifford}

Soient $a,b \in \R^*_+$ et $\T_{a,b} \subset \C P^2$ le tore lagrangien d\'efini par les \'equations  $|x| = a$, $|y|=b$
dans les coordonn\'ees affines $(x,y) \in (\C)^2 \subset \C P^2$. 
Ce tore est le lieu fixe de l'involution antibirationnelle de Cremona 
$c^{a,b} : (x, y , z) \in \C P^2 \setminus \{ (1, 0 , 0) , (0 , 1 , 0) , (0, 0 , 1) \} \mapsto (a^2 \overline{yz}, b^2 \overline{xz} ,  \overline{xy}) \in \C P^2$.
Cette involution antibirationnelle $c^{a,b}$, $a,b \in \R^*_+$, est simple. 
En effet, soit $Y$ le plan projectif \'eclat\'e aux trois points 
$(1, 0 , 0) , (0 , 1 , 0) , (0, 0 , 1)$ et $\pi : Y \to \C P^2$ la projection associ\'ee. L'involution $c^{a,b}$ se rel\`eve en une involution antiholomorphe
$c^{a,b}_Y$ d\'efinie partout, soit une structure r\'eelle. 
Cette derni\`ere envoie les trois diviseurs exceptionnels sur les transform\'ees strictes des c\^ot\'es du  triangle $(1, 0 , 0) , (0 , 1 , 0) , (0, 0 , 1)$, d'o\`u la simplicit\'e de $c^{a,b}$. 
Ainsi, le Th\'eor\`eme \ref{theoantibir} s'applique et fournit des invariants $\chi_r^{d_Y}$ par d\'eformation du triplet  $(\C P^2 , \omega , c^{a,b})$. 
Le deuxi\`eme groupe d'homologie de $Y$ est engendr\'e par une droite g\'en\'erique et les diviseurs exceptionnels $E_1, \dots E_3$ de nos \'eclatements.  
La classe d'homologie $d_Y$ de nos courbes rationnelles r\'eelles de $Y$ est d\'etermin\'ee par quatre entiers $d, d_1 , \dots , d_3$ satisfaisant
la relation $d = d_1 + d_2 + d_3$. Si l'on contracte $E_1, \dots E_3$, ces courbes se contractent sur des courbes rationnelles de degr\'e $d$ du plan
qui ont un point de multiplicit\'e $d_1 , d_2 , d_3$ en $(1, 0 , 0) , (0 , 1 , 0) , (0, 0 , 1)$ respectivement. 
Ces courbes rationnelles immerg\'ees ont en outre la propri\'et\'e de rencontrer le tore  $\T_{a,b}$ 
en une collection de points isol\'es et en un cercle immerg\'e, elles consistent en fait en une paire de disques $J$-holomorphes qui reposent sur $\T_{a,b}$
et sont \'echang\'es par $c^{a,b}$. Si l'on contracte plut\^ot un diviseur exceptionnel, disons $E_3$, ainsi que son image sous $c^{a,b}_Y$, 
alors on obtient des courbes rationnelles de bidegr\'e  $(d_1 + d_3 , d_2 + d_3)$ sur l'hyperbolo\"{\i}de quadrique $(\C P^1 \times \C P^1 , \text{conj} \times \text{conj})$,
qui ont une paire de points de multiplicit\'e $d_3$ en deux points complexes conjugu\'es, \`a savoir les points o\`u $E_3$ et son image se contractent. Lorsque $d_3 = 0$ ou $1$, 
cet invariant $\chi_r^{d_Y}$ vaut l'invariant correspondant dans l'hyperbolo\"{\i}de quadrique $(\C P^1 \times \C P^1 , \text{conj} \times \text{conj})$, \`a savoir 
$\chi_r^{(d_1, d_2)}$ et $\chi_r^{(d_1 + 1, d_2 + 1)}$
respectivement. Des estimations de ces derniers se trouvent dans \cite{IKS}.

\begin{cor}
Soient $r , s , d \in \N$ tels que $ r+ 2s = 2d - 1$ et supposons donn\'ee une collection de $r$ points distincts dans $\T_{a,b} \subset \C P^2$, $a,b \in \R^*_+$ 
ainsi qu'une collection de $s$
paires distinctes de points dans $ \C P^2 \setminus \big( \T_{a,b} \cup \{ (1, 0 , 0) , (0 , 1 , 0) , (0, 0 , 1) \} \big)$ \'echang\'ees par l'involution antibirationnelle $c^{a,b}$.
Alors, pour tous $d_1 , d_2 , d_3 \in \N$ tels que $d = d_1 + d_2 + d_3$, il y a au moins $\vert \chi_r^{d_Y} \vert$  paires de disques $J_X$-holomorphes reposant sur $\T_{a,b}$,
\'echang\'es par $c^{a,b}$, passant par les  $r$ points donn\'es et intersectant chacune des $s$ paires de points complexes conjugu\'es, d\`es lors que $J_X$
se rel\`eve en une structure $J_Y$ appartenant \`a l'un des $\R {\cal J}_{\omega_Y}^{U , J_0}$ donn\'e par le Lemme \ref{lemmaJ0}. La r\'eunion de ces deux disques
dans chacune de ces paires forme une courbe rationnelle plane de degr\'e $d$ ayant un point de multiplicit\'e $d_1 , d_2 , d_3$ en $ (1, 0 , 0) , (0 , 1 , 0) , (0, 0 , 1)$ 
respectivement. $\square$
\end{cor}

Remarquons que des invariants \'enum\'eratifs portant sur des disques \`a bords dans le tore de Clifford ont \'et\'e obtenus par
 P. Biran et O. Cornea \cite{Biran}. Les disques holomorphes \`a bords dans le tores de Clifford ont par ailleurs \'et\'e \'etudi\'es par C.-H. Cho dans sa th\`ese en termes de produits de  Blaschke. 
En ce qui concerne nos r\'esultats pr\'esent\'es dans ce paragraphe \ref{subsectpre}, il reste \`a s'affranchir de la notion de simplicit\'e (des involutions antibirationnelles).

\vspace{0.7cm}
\noindent
Universit\'e de Lyon ; CNRS ;
Universit\'e Lyon~1 ; Institut Camille Jordan

\end{document}